\documentclass{article}
\usepackage[utf8]{inputenc}
\usepackage[english]{babel}
\usepackage{appendix}
\usepackage{color}
\usepackage{amsmath,amssymb,amsthm,yhmath,amsfonts,mathrsfs,enumitem,mathtools,dsfont}
\usepackage{multirow}
\usepackage{subcaption}
\usepackage{multirow}

\newtheorem{theorem}{Theorem}

\renewcommand{\d}{\mathrm{d}}

\title{Mean-field analysis of load balancing principles in large scale systems}

\author{
	Ill\'es Horv\'ath\\
	{\small ELKH-BME Information Systems Research Group}\\
	{\small e-mail: horvath.illes.antal@gmail.com}\\[3mm]
	M\'arton M\'esz\'aros \\
	{\small Department of Stochastics, Mathematics Institute}\\
	{\small Budapest University of Technology and Economics}\\
	{\small e-mail: mezsimarci@gmail.com}
}

\date{}

\begin{document}
	
	\maketitle
	
	\begin{abstract}
		
		Load balancing plays a crucial role in many large scale systems. Several different load balancing principles have been proposed in the literature, such as Join-Shortest-Queue (JSQ) and its variations, or Join-Below-Threshold. We provide a high level mathematical framework to examine heterogeneous server clusters in the mean-field limit as the system load and the number of servers scale proportionally. We aim to identify both the transient mean-field limit and the stationary mean-field limit for various choices of load balancing principles, compute relevant performance measures such as the distribution and mean of the system time of jobs, and conduct a comparison from a performance point of view.
		
	\end{abstract}
	
\section{Introduction}

For large scale service systems, where service resources (e.g. computing capacity) are distributed to several service units, load balancing plays a crucial role in distributing the total load of the system to ensure better overall service for the incoming tasks (jobs).

There are many different types of load balancing principles. Static load balancing does not take into account the state of the system, instead aiming for a balanced distribution based purely on the incoming jobs. Static load balancing is in general easy to set up, requires minimal overhead communication and performs well when the incoming jobs have some regular patterns.

However, in most systems the incoming jobs have some level of random variability. This situation is generally better handled by load balancing policies which take into account the current state of the system. Scheduling decisions may be based on different types of information, depending on what is available. In general, one of the most important parameters is the current load of the servers, as it is generally desirable to maintain a balanced load among all servers. If available, further information taken into account may include any of the following:
\begin{itemize}
	\item the servers may be heterogeneous, with faster and slower servers;
	\item job and server types may be important in case the servers are heterogeneous and certain servers can serve certain types of jobs more efficiently;
	\item job sizes may be used to compute current server load more precisely;
	\item in some cases, physical location may play a role;
	\item there may be bottlenecks other than computing capacity in the system (e.g. bandwidth).
\end{itemize}

In many real-life systems, such information may not be available, but even if it is, there is a tradeoff: a complicated load balancing policy that requires too much communication and computation may generate a significant overhead cost, slowing down the entire system. Hence it is in general desirable to stick to simple load balancing policies. In the present paper, we provide a mathematical framework that does not include communication overhead costs. Such aspects can be addressed in the modeling in several ways; however, these are highly scenario-dependent, and as such, we decided to keep the model high-level.

We will discuss load balancing policies based exclusively on the queue length of servers. Job types, physical location and other bottlenecks will not play a role. We allow a heterogeneous server cluster, where there are several different types of servers, and the model can also incorporate processor sharing, where a server can serve multiple jobs simultaneously.

The server cluster model of the present paper will be described by a \emph{density-dependent Markov population model}. As the system size goes to infinity, the mean-field limit of density-dependent Markov population models has been examined in the literature for both the transient regime (up to a finite time horizon) and in the stationary regime.

The transient limit object is deterministic and can be described as the solution a system of ordinary differential equations (ODEs) in case the Markov transition rates are Lipschitz-continuous \cite{Kurtz1}, or as the solution of a differential inclusion in case the transition rates are discontinuous \cite{Gast1}. Overall, these results are relatively straightforward to apply for the model in the present paper.

For the stationary regime, for Lipschitz-continuous transition rates, it is known that in the mean-field limit, the stationary distribution of the finite system concentrates on the unique asymptotically stable solution (attractor) of the limit system of ODEs \cite{Kurtz2}. Similar results available for the discontinuous setting, but only in case the attractor lies inside a domain where the transition rates are continuous \cite{Gast1}. We are not aware of any general results in case the attractor is at a discontinuity point of the transition rates, which happens to be the case for several of the load balancing policies discussed in the present paper.

The contributions of the paper are the following:
\begin{enumerate}
	\item Providing a high-level mathematical framework for modelling load balancing systems that accommodates several different load balancing principles.
	\item Identification of the mean-field limit in both the transient and stationary regime. 
    \item Computation of the mean service time and also the service time distribution in the stationary mean-field limit. Computation techniques need to be adapted for discontinuities; these modified formulas are, to the best of our knowledge, novel.
	\item Numerical comparison of the various load balancing principles via simulation and theoretical computations for the mean-field limit.
\end{enumerate}

All of the above is carried out for a fairly general setting, where the server cluster can be heterogeneous, and we will also allow a varying service rate, depending on the number of jobs in a given server. We will focus mostly on first-in-first-out (FIFO) service principle, but note that all calculations are straightforward to derive for limited processor sharing (LPS), where a server can serve multiple jobs simultaneously. 

Rigorous proofs are not the main focus of the paper. We do refer to relevant rigorous results from the literature in cases where they are available, but only provide heuristic arguments for the novel cases. That said, numerical analysis does support the heuristic computations of the paper.

The codes used for the simulations and analytic calculations throughout the paper are available at \cite{code}.

The rest of the paper is structured as follows: the rest of this section is dedicated to an overview of load balancing in the literature (Section \ref{ss:lbp}), and to the necessary mathematical background in queueing theory (Section \ref{ss:bdp}) and population processes (Section \ref{ss:pop}). Section \ref{s:cluster} describes the general setup of the server cluster we are interested in. Section \ref{s:loadbalance} describes the various load balancing principles. Section \ref{s:numexp} contains numerical experiments and comparison of the various load balancing principles, and Section \ref{s:concl} concludes the work. The Appendix addresses a few related questions not strictly part of the main body of work, and also some further details.

\subsection{Load balancing principles}
\label{ss:lbp}

One of the classic dynamic load balancing policies is Join-Shortest-Queue (JSQ), where the incoming job is assigned to the server with the shortest queue (lowest number of jobs) \cite{jsq}. The upside of this method is that it offers very even balancing for homogeneous server clusters. However, it requires up-to-date knowledge of all server states, which may require a significant communication overhead.

Due to this, several variants of JSQ have been in use: for JSQ($d$), the incoming job is scheduled to the shortest queue from among $d$ servers, selected at random. This offers less balanced load distribution, but also requires less communication. $d=1$ corresponds to random assignment with no load balancing, and $d$ equal to the total number of servers corresponds to JSQ; as $d$ is increased, it offers better balancing but also more overhead communication. Interestingly, already for $d=2$, the resulting load balancing policy has certain asymptotic optimality properties \cite{mitzenmacher2}, often referred to as the power-of-2 (or power-of-$d$) policies. As a consequence, $d$ is often selected relatively low, such as $d=2$ or $d=5$. 

For Join-Idle-Queue (JIQ), the incoming job is scheduled to an idle server at random; if there are no idle servers, the assignment is random among all servers. Once again, this offers less balanced load distribution and less communication overhead than JSQ, but, similar to JSQ($d$), has some nice asymptotic optimality properties. Mean-field analysis has been carried out for JIQ in \cite{Mitzenmacher}.

Another related load balancing policy is Join-Below-Threshold (JBT), which associates a threshold with each server; servers below their threshold are considered available and servers at or above their threshold are full. Jobs will be dispatched to a server randomly from among all available servers. This policy again offers less balancing than JSQ, but still offers protection against overloaded servers, and requires communication only when a server switches between available and full. For a full mean-field analysis and cluster optimization of JBT, we refer to \cite{jbt}.

\subsection{Birth-death processes and queues}
\label{ss:bdp}

The jobs arriving to and leaving a server's queue can be modelled with a birth-death process (Markov-queue).

For technical simplicity, we resort to finite queues, with the maximal queue length denoted by $B$ and state space of a single queue $\Omega =  \{0,1,2,\dots,B\}$.

We assume Markov arrivals, that is, jobs arrive according to a Poisson process, and Markov service, that is, the time it takes to serve a job (once service has started) is exponentially distributed. 

There are multiple service principles. For First-In-First-Out (FIFO) service principle, the server always serves the first job of a queue, while the other jobs wait. Whenever the first job has finished service, the server immediately starts serving the next job in the queue. For Limited Processor Sharing (LPS), the server can work on multiple jobs simultaneously. The maximum number of jobs served simultaneously is called the multi-programming level (MPL); further jobs in the queue wait and enter service in a manner similar to FIFO. We allow the service rate to depend on the number of jobs in the queue (this is particularly relevant for LPS, where multiple jobs can be served jointly for more efficient service overall). The choice of service principle has no effect on the queue length changes (no matter which job is served, queue length decreases by 1), but it does affect the system time of individual jobs. We will mostly focus on FIFO.

\subsection{Density-dependent population processes}
\label{ss:pop}

In this section, we present mathematical background and framework for density-dependent Markov population processes.


A density-dependent Markov population process has $N$ interacting components, each of which is in a state from a finite set of local states $S$. The global state of the system is defined as the total number of individuals in each state, that is,
a vector
$X^N\in \{0,1,\dots,N\}^{|S|}$
with
$X^N_1+\dots+X^N_{|S|}=N$. The normalized global
state of the system can be defined as
$$x^N=\frac{X^N}{N},$$
so
$x^N\in [0,1]^S$ with $x^N_1+\dots+x^N_{|S|}=1.$

Each component acts as a continuous time Markov chain. The rate of the transition from $i \in S$ to $j \in S$ is $r_{ij}^N$ (for $i \neq j$). The rates are assumed to be \emph{density-dependent}, that is $$r_{ij}^N = r_{ij}(x^N)$$ for some function $r_{ij}:[0,1]^{|S|}\to[0,\infty]$. In the classic setup defined by Kurtz \cite{Kurtz1,Kurtz2}, the functions $r_{ij}$ are usually assumed to be Lipschitz-continuous and independent of $N$. With this setup, $x^N(t)$ is a continuous time Markov-chain. We define the \emph{mean-field equation} of the system as the following:
\begin{align}
\label{eq:mf0}
\frac{\d}{\d t}v_i(t)=\sum_{j\in S} v_j(t)r_{ji}(v(t)),\quad i\in S,
\end{align}
where $$r_{ii}:=-\sum_{j\in S, j\neq i}r_{ij},$$
and $$x_i^N(0)\to v_i(0)\quad (\textrm{for }i=1,\dots, |S|),\quad\textrm{in probability as }N\to\infty.$$
Lipschitz-continuity guarantees existence and uniqueness of the solution of \eqref{eq:mf0}. The following result of Kurtz states mean-field convergence in the transient regime \cite{Kurtz1,Kurtz2,Kurtz3}:
\begin{theorem}[Transient mean-field convergence]
	\label{t:kurtz-trans}
	Assuming $r_{ij}$ $ (i,j\in S)$, are Lipschitz-continuous and 
	$$x_i^N(0)\to v_i(0)\quad i\in \{1,\dots,|S|\} ,\quad\textrm{in probability},$$
	then for any $T>0$ we have 
	\begin{equation*}
	\lim_{N \rightarrow \infty}
	\mathds{P}\left(
	\sup_{t \in [0,T]} \|\bar{\mathbf{x}}^N(t) - \mathbf{v}(t)\| > \epsilon\right) = 0.
	\end{equation*}
\end{theorem}
Kurtz also proved that the standard deviation of $x^N$ is of order $\frac{1}{\sqrt{N}}$ \cite{Kurtz2}.

An important concept related to Theorem \ref{t:kurtz-trans} is asymptotic independence, also known as propagation of chaos, stating that as $N\to\infty$, the evolution of two distinct queues is asymptotically independent. This is due to the fact that the evolution of a queue depends only on the global state, which is asymptotically deterministic.

We also have stationary mean-field convergence.
\begin{theorem}[Stationary mean-field convergence]
	\label{t:kurtz-stat}
	Given the following assumptions:
	\begin{itemize}
		\item $r_{ij}$ are Lipschitz-continuous,
		\item the Markov process $x^N(t)$ has a unique stationary distribution $\pi^N$
		for each $N$, and
		\item \eqref{eq:mf0} has a unique stable attractor $(\nu_1,\dots,\nu_{|S|})$,
	\end{itemize}
	we have that the probability measure $\pi^N$ on $S$ converges in probability to the Dirac measure concentrated on $\nu$.
\end{theorem}

Theorems \ref{t:kurtz-trans} and \ref{t:kurtz-stat} have been generalized in several directions during recent years. Bena\"{i}m and Le Boudec elaborated a framework applicable for a wider range of stochastic processes, which also allows the $r_{ij}$ functions to have a mild dependency on $N$ \cite{Benait1}.

The condition on Lipschitz-continuity can also be weakened. For discontinuous $r_{ij}$'s, \eqref{eq:mf0} turns into a differential inclusion. A formal setup for differential inclusions is quite technical, and is omitted from the present paper. For a fully detailed setup, we refer to \cite{Gast1}, specifically Theorems 4 and 5, and \cite{roth2013stochastic}, Theorem 3.5 and Corollary 3.9 for a corresponding version of Theorem \ref{t:kurtz-trans}.

For a corresponding version of Theorem \ref{t:kurtz-stat} for discontinuous transition rates, we refer to \cite{Gast1}, where the main additional condition is that the unique attractor lies inside a domain where the $r_{ij}$ are continuous.

The applicability of Theorems \ref{t:kurtz-trans} and \ref{t:kurtz-stat} will be addressed more in Section \ref{s:cluster}.

From Theorem \ref{t:kurtz-stat} it also follows that $$\lim_{N\to\infty} \mathds{E}(\pi^N)= \nu,$$
so $\nu$ can be used as an approximation for $\mathds{E}(\pi^N)$ for large $N$. $\mathds{E}(\pi^N)$ here is basically an $|S|$-dimensional vector of distributions, which converges to a constant $|S|$-dimensional vector in distribution. The limit point can be interpreted as a distribution on $S$, and is the stable attractor $\nu$.

\section{Server clusters}
\label{s:cluster}

The server cluster model examined in the present paper consists of $N$ servers, each with a finite buffer, and a single common dispatcher. Jobs arrive to the dispatcher according to a Poisson process with rate $N\lambda$ (that is, the average arrival rate is $\lambda$ per server). Each arriving job is instantly dispatched to one of the $N$ servers; that is, the dispatcher maintains no queue.

The cluster may have $K$ different \emph{server types}. We assume $K$ is fixed, independent from $N$.

The servers within each type are identical. Buffer sizes are denoted by $B^{(k)}$ for each type $k\in \{1,\dots,K\}$. We assume service times are exponentially distributed; for each server type, the service rate can be constant or it may depend on the current queue length of the server.
Service rates are denoted by $\mu_i^{(k)}$, where $i\in\{ 0,1,\dots,B^{(k)}\}$ is the queue length, and $k\in\{ 1,2,\dots,K\} $ denotes the type of the server. For a given $k\in \{1,\dots,K\}$, $\mu_0^{(k)},\dots,\mu_{B^{(k)}}^{(k)}$ is also referred to as the \emph{service rate curve}. ($\mu_0^{(k)}=0$, but we still include it in the notation.)

For each service rate curve, it is natural to assume that the total rate increases with the queue length, but the per-job rate decreases with the queue length:
\begin{align}
\label{eq:mon}
\mu_1^{(k)}\leq \mu_2^{(k)}\leq\mu_3^{(k)}\leq\dots, \quad \mu_1^{(k)}\geq \frac{\mu_2^{(k)}}2\geq\frac{\mu_3^{(k)}}3\geq\dots \quad k\in\{ 1,2,\dots,K\}
\end{align}

Due to the finite buffer sizes, data loss may occur whenever a job is dispatched to a full queue. The probability of a job loss will be typically very low (due to load balancing), but it is still something that we will address in due course.

The server cluster is a density-dependent population process, where the state of a server is simply the number of jobs in its queue. The global state will be denoted by
$$X_{i}^{(k),N}(t),\qquad \left( 0\leq i\leq B^{(k)},\, 1\leq k\leq K \right),$$
where $X_{i}^{(k),N}(t)$ is the number of servers with $i$ jobs in its queue at time $t$. We will mostly use its normalized version $$x^{N}(t)=x_{i}^{(k),N}(t),\qquad \left( 0\leq i\leq B^{(k)},\qquad 1\leq k\leq K\right),$$ where $$x_{i}^{(k),N}(t)=\frac{X_{i}^{(k),N}(t)}{N}.$$
The number of servers of type $k$ is denoted by $N_k$ and the ratio of each server type is denoted by 
$$\gamma_k^N=\frac{N_k}N, \, \qquad k=1,\dots,K.$$
$\gamma_k^N$ may depend on $N$, but we will assume they converge to some fixed values $\gamma_k$ as $N\to\infty$. We also want the system to be stable, so 
\begin{align}
\label{eq:stable}
\lambda < \sum_{k=1}^K \gamma_k^N \mu_B^{(k)}.
\end{align}
(Actually, due to the finite buffer size assumption, the system is technically always stable, but we will nevertheless assume \eqref{eq:stable}.)

The evolution of $x^N(t)$ can be formally defined using Poisson representation. Let 
\begin{align*}
P_{i\to (i+1),k}(t),&\quad 0\leq i\leq B^{(k)}-1,\, k=1,\dots,K\\
P_{i\to (i-1),k}(t),&\quad 1\leq i\leq B^{(k)},\,k=1,\dots,K
\end{align*}
denote independent Poisson processes with rate 1. $P_{i\to (i+1),k}(t)$ corresponds to arrivals to queues of type $k$ with length $i$, and $P_{i\to (i-1)}(t)$ corresponds to jobs leaving queues of type $k$ with length $i$.

The Poisson representation of $x^{N}(t)$ is
\begin{align}
\label{eq:poirepr}
\begin{split}
x_{i}^{(k),N}(t)=&\frac{1}{N} P_{(i-1)\to i,k}\left(N \int_0^t \lambda f^{(k)}_{i-1}(x^N(s))\d s\right)\\
&\quad -\frac{1}{N} P_{i\to (i+1),k}\left(N \int_0^t \lambda f^{(k)}_i(x^N(s))\d s\right)\\
&\quad +\frac{1}{N} P_{(i+1)\to i,k}\left(N \int_0^t \mu^{(k)}_{i+1} x_{i+1}^{(k),N}(s)\d s\right)\\
&\quad -\frac{1}{N} P_{i\to (i-1),k}\left(N \int_0^t \mu^{(k)}_i x_{i}^{(k),N}(s)\d s\right),
\end{split}
\end{align}
where $f_i^{(k)}(x^N(t))$ is the probability of a new arriving job to enter a queue with length $i$ of type $k$. 

The
$$\{f_i^{(k)}(x^N(t)):0\leq i\leq B_k,\, k=1,\dots,K\}$$
functions are going to be collectively called the \emph{dispatch functions}. The dispatch functions depend on the load-balancing principle, which will be addressed later. Formally, $f_i^{(k)}$ are defined on the normalized state $x^{N}(t)$, which are all contained in the domain
\begin{align}
\label{eq:domain}
\{x:x\in\mathbb{R}^{\sum_{k=1}^K (B^{(k)}+1)},\, x_j^{(k)}\geq 0,\,\sum_{k=1}^K \sum_{j=0}^{B^{(k)}}x_j^{(k)}=1\}.
\end{align}

The four possible changes in the number of queues of length $i$ which appear in (\ref{eq:poirepr}) correspond to:
\begin{itemize}
\item a job arriving to a queue of length $i-1$;
\item a job arriving to a queue of length $i$;
\item a job leaving a queue of length $i+1$;
\item a job leaving a queue of length $i$.
\end{itemize}
On the border of the domain \eqref{eq:domain}, certain changes cannot occur. There is no service in empty queues: $$\mu_0^{(k)}=0\qquad (k=1,\dots,K),$$
and no arrival to full queues:
$$f^{(k)}_{B^{(k)}}(.)\equiv 0\qquad (k=1,\dots,K).$$

We are interested in server clusters of various $N$ sizes and especially the limit object as $N\to\infty$, that is, the mean-field limit (in accordance with Section \ref{ss:pop}). We first define the general mean-field equations corresponding to (\ref{eq:poirepr}):
\begin{align}
\label{eq:mf}
\begin{split}
v^{(k)}_{i}(t)=\,\,&v^{(k)}_i(0)+\int_0^t \lambda f^{(k)}_{i-1}(v(s))\d s -\int_0^t \lambda f^{(k)}_i(v(s))\d s\\ &\quad +\int_0^t \mu^{(k)}_{i+1} v_{i+1}^{(k)}(s)\d s -\int_0^t \mu^{(k)}_i v_{i}^{(k)}(s)\d s
\end{split}
\end{align}
in integral form, or, equivalently,
\begin{align}
\label{eq:mf_diff}
\frac{\d}{\d t}v_{i}^{(k)}(t)=\lambda f^{(k)}_{i-1}(v(t))-\lambda f_i^{(k)}(v(t))
+\mu^{(k)}_{i+1} v^{(k)}_{i+1}(t) - \mu^{(k)}_i v^{(k)
}_{i}(t)
\end{align}
in differential form. An empty initial cluster corresponds to the initial condition
\begin{align*}
v_i^{(k)}(0)=
\left\{
\begin{array}{l}
\gamma_k\quad\, \textrm{ for }i=0,\\
0\qquad \textrm{otherwise.}
\end{array}
\right.
\end{align*}

Theorem \ref{t:kurtz-trans} applies to this system whenever the $f_i^{(k)}$ functions are Lipschitz-continuous. It turns out that the conditions of the general version of Theorem \ref{t:kurtz-trans} are mild enough so that transient mean-field convergence holds for all the discontinuous choices of $f_i^{(k)}$ in the present paper, but this is not checked rigorously.

For the stationary case, we denote the stationary distribution 
$$\nu=(\nu_i^{(k)}),i=0,\dots,B^{(k)},\quad k=1,\dots,K$$
(similar to the notation of Section \ref{ss:pop}). Theorem \ref{t:kurtz-stat} applies whenever $f_i^{(k)}$ are Lipschitz-continuous. In the discontinuous setting, the most relevant question is whether the $f_i^{(k)}$ functions are continuous at the unique fixed point $\nu$ or not. If $\nu$ lies inside a region where $f^{(k)}_i$ are Lipschitz-continuous, then the conclusion of Theorem \ref{t:kurtz-stat} applies. However, when  the $f_i^{(k)}$ functions are discontinuous at $\nu$, Theorem \ref{t:kurtz-stat} does not apply; in fact, little is known in this case rigorously. Based on this, it makes sense to distinguish the following two cases:
\begin{enumerate}
	\item the functions $f_i^{(k)}$ are Lipschitz-continuous at $\nu$, or
	\item the functions $f_i^{(k)}$ are discontinuous at $\nu$.
\end{enumerate}

When the functions $f_i^{(k)}$ are Lipschitz-continuous at $\nu$, the equations for the mean-field stationary distribution can be obtained from \eqref{eq:mf_diff} by setting $\frac{\d}{\d t}v^{(k)}_i(t)=0$:
\begin{align}
\label{eq:mf_stat}
\begin{split}
0=\lambda f^{(k)}_{i-1}(v(t))-\lambda f^{(k)}_i(v(t))
+\mu^{(k)}_{i+1} v^{(k)}_{i+1}(t) - \mu^{(k)}_i v^{(k)}_{i}(t)\\
i\in \{1,\dots,B^{(k)-1} \} ,\,\qquad k\in\{ 1,\dots,K \}
\end{split}
\end{align}
which are equivalent to the dynamic balance equations
\begin{align}
\label{eq:dynamicbalance2}
\mu^{(k)}_i \nu^{(k)}_i =\lambda f^{(k)}_{i-1}(\nu),\qquad i\in \{1,\dots,B^{(k)}\} ,\, \qquad k\in \{ 1,\dots,K \} .
\end{align}
We also have equations for the ratio of each server type:
\begin{align}
\label{eq:servertypes}
\sum_{i=0}^{B^{(k)}} \nu^{(k)}_i=\gamma_k, \qquad k\in \{1,\dots,K \}.
\end{align}
\eqref{eq:dynamicbalance2} + \eqref{eq:servertypes} provide algebraic equations for $\nu$.

We also propose another approach to obtain $\nu$ numerically, by solving the transient equations \eqref{eq:mf_diff} and taking the solution at a large enough point in time. (This assumes convergence to a single asymptotically stable solution, which we do not aim to prove rigorously.)

When the $f^{(k)}_i$ are discontinuous at $\nu$, more considerations are needed to derive the dynamic balance equations. This will be addressed separately for each load balancing principle.

Further remarks.

The assumption that both arrival and service are Markovian means that the entire system is a Markov (population) process, which keeps the setup fairly simple. Interestingly, the same mean-field limit would be obtained for any arrival process as long as the arrivals average out in the mean-field limit; to be more precise, for any arrival process for which the Functional Strong Law of Large Numbers holds (see e.g. Theorem 3.2.1 in \cite{whitt_suppl}).

In case the monotonicity condition \eqref{eq:mon} does not hold, mean-field convergence may fail. \cite{jbt} contains specific examples where \eqref{eq:mf} has multiple fixed points; stable fixed points correspond to quasi-stationary distributions of the population process for any finite $N$. The solution of \eqref{eq:mf} will converge to one of the stable fixed points (depending on the initial condition). However, for any finite $N$, the population process will spend very long periods of time near one of the quasi-stationary points, switching between these points infinitely often.

\subsection{Mean system time}
A wide variety of parameters can be considered to describe the efficiency of such a system. A natural choice is the mean system time: the average time a job spends in the system between its arrival and service. We aim to calculate the mean system time $H$ in the stationary mean-field regime. We note that the mean system time is a somewhat artificial object here since technically there are no individual jobs in the mean-field limit. It may be helpful to think of the mean-field limit as the case when $N$ is extremely large.

One way to compute $H$ is via Little's Law
\begin{align*}
    H=L/\lambda_e,
\end{align*}
where $L$ is the mean queue length in the system, and $\lambda_e$ is the effective arrival rate (which excludes jobs not entering the system due to job loss). From the mean-field stationary distribution $\nu$, $L$ is easily computed, while $\lambda_e$ depends on the load balancing policy, but is typically also straightforward to compute. Little's law can actually be applied to each server type separately for more detailed information; this is addressed in Appendix \ref{app:little}.

Here we propose a different method to compute $H$, which gives even more detailed information, and will be useful later on. Let $H_{i,j}^{(k)}$ denote the mean time until service for a job that is in position $i$ in a queue of type $k$ with $j$ jobs total (so $1\leq i \leq j \leq B^{(k)}$, $1\leq k\leq K$).

In the case of constant service rates, $H^{(k)}_{i,j}= \frac{i} {\mu^{(k)}}$ holds. For non-constant service rate curves however, the service rate may change due to later arrivals, so we need to keep track of both the length of the queue and the position of the job within it. We will derive a system of linear equations using total expectation and the Markov property. For simplicity, we assume FIFO service principle in the following calculations, but due to Little's law, this assumption does not affect the value of $H$.

The mindset is that we are following a tagged job at position $i$ of a queue of type $k$ with total queue length $j$, and the equations are based on possible changes in the queue, with the environment fixed due to the stationary mean-field regime.

\begin{align}
\nonumber
H_{i,j}^{(k)} &= \frac{1}{\lambda f_{j}^{(k)}(\nu)/{\nu_j^{(k)}}+\mu_{j}^{(k)}}+
\frac{\lambda {f_{j}^{(k)}(\nu)}/{\nu_j^{(k)}}}{\lambda {f_{j}^{(k)}(\nu)}/{\nu_j^{(k)}}+\mu_{j}^{(k)}}H_{i,j+1}^{(k)}+\\
\nonumber
& \qquad \quad
\frac{\mu_{j}^{(k)}}{\lambda {f_{j}^{(k)}(\nu)}/{\nu_j^{(k)}}+\mu_{j}^{(k)}}H_{i-1,j-1}^{(k)}
\qquad  (2\leq i\leq j\leq B^{(k)}-1),\\
\label{eq:Hij}
H_{i,B^{(k)}}^{(k)} &=\frac{1}{\mu_{B^{(k)}}^{(k)}}+H_{i-1,B^{(k)}-1}^{(k)} \qquad (2\leq i\leq B^{(k)}),\\
\nonumber
H_{1,j}^{(k)} &=\frac{1}{\lambda {f_{j}^{(k)}(\nu)}/{\nu_j^{(k)}}+\mu_{j}^{(k)}}+\frac{\lambda {f_{j}^{(k)}(\nu)}/{\nu_j^{(k)}}}{\lambda {f_{j}^{(k)}(\nu)}/{\nu_j^{(k)}}+\mu_{j}^{(k)}}H_{1,j+1}^{(k)} \\ \nonumber & \qquad \qquad \qquad (1\leq j\leq B^{(k)}-1),\\
\nonumber
H_{1,B^{(k)}}^{(k)} &=\frac{1}{\mu_{B^{(k)}}^{(k)}}.
\end{align}
\eqref{eq:Hij} makes use of the standard one step argument. We focus on a single queue of a given type $k$ in the mean-field limit while assuming the environment to be stationary, and look for the next possible change in that queue. Jobs arrive to type $k$ servers of queue length $j$ with a rate of $N\lambda f_j^{(k)}(\nu)$, and each job will be sent to one of $N\nu_j^{(k)}$ servers, so the arrival rate at a specific queue will be
$$\frac{N\lambda f_j^{(k)}(\nu)}{N\nu_j^{(k)}}=\lambda{f_{j}^{(k)}(\nu)}/{\nu_j^{(k)}},$$
while the service rate is $\mu_{j}^{(k)}$, so the rate of any change for a queue of length $j$ is $\lambda {f_{j}^{(k)}(\nu)}/{\nu_j^{(k)}}+\mu_{j}^{(k)}$. The change will either increase or decrease the length of the queue by 1, and we can apply total expectation.

For full queues ($j=B^{(k)}$), arrival is not possible, that is, $f_{B^{(k)}}^{(k)}(.)\equiv 0$ for $k=1,\dots,K$.

In order to solve \eqref{eq:Hij}, we first obtain the mean-field stationary distribution $\nu$. $\nu$ can be calculated from either the balance equations \eqref{eq:mf_stat} when possible, or by numerically solving the transient mean-field equations \eqref{eq:mf_diff} and setting $t$ large enough. Once $\nu$ is obtained, \eqref{eq:Hij} is just a system of linear equations for $H_{i,j}^{(k)}$, which can actually be solved separately for each $k$ for $1\leq k \leq K$. Once \eqref{eq:Hij} is solved, the mean system time $H$ is just a linear combination of the values $H_{j,j}^{(k)}$ according to the probabilities with which a job will be scheduled to a queue of length $j-1$ of a $k$-type server, that is,
\begin{align}
\label{eq:Hweights}
H=\frac{1}{\sum_{k=1}^K\sum_{j=1}^{B^{(k)}} f^{(k)}_{j-1}(\nu)} \sum_{k=1}^K\sum_{j=1}^{B^{(k)}} f^{(k)}_{j-1}(\nu)H_{j,j}^{(k)}.
\end{align}

The normalizing factor in \eqref{eq:Hweights} addresses job loss, as we only want to consider the mean system time of jobs which actually enter the system. Job loss probability is equal to
$$1-\sum_{k=1}^K\sum_{j=1}^{B^{(k)}} f^{(k)}_{j-1}(\nu).$$

\eqref{eq:Hij} and \eqref{eq:Hweights} are only valid if the dispatch functions $f_i^{(k)}$ are continuous at $\nu$. In other cases, we may need to tweak the formulas. We will provide the corresponding versions of \eqref{eq:Hij} and \eqref{eq:Hweights} on a case-by-case basis whenever the functions $f_i^{(k)}$ are discontinuous at $\nu$. These versions will be heuristic in the sense that no formal rigorous proof will be provided, but the results nevertheless agree with the results from simulations.

\subsection{System time distribution}
\label{ss:systemtime}

In this section, we calculate the system time distribution for a random job. Here, the service principle is actually important; we will present the calculation for FIFO service principle here. The calculations need to be modified for LPS service principle; the corresponding equations are provided in Appendix \ref{app:lps}.

Let $h_{i,j}^{(k)}(t)$ denote the probability density function of the remaining system time of a job at position $i$ in a queue of length $j$ and queue type $k$. Its Laplace-transform is defined as
$$\tilde H_{i,j}^{(k)}(s)=\int_0^\infty h_{i,j}^{(k)}(t)e^{-st}\mathrm{d}t.$$

The following system of equations is the corresponding version of \eqref{eq:Hij} for the Laplace-transforms instead of the means. Total expectation also applies to Laplace-transforms, and we use the fact that the Laplace-transform of 0 is 1 and the Laplace-transform of $\lambda e^{-\lambda t}$ is $\frac{\lambda}{s+\lambda}$ to obtain
\begin{align}
\nonumber
\tilde H_{i,j}^{(k)}(s) &= \frac{\lambda f_{j}^{(k)}(\nu)/{\nu_j^{(k)}}+\mu_{j}^{(k)}}{s+\lambda f_{j}^{(k)}(\nu)/{\nu_j^{(k)}}+\mu_{j}^{(k)}}\Bigg(
\frac{\lambda {f_{j}^{(k)}(\nu)}/{\nu_j^{(k)}}}{\lambda {f_{j}^{(k)}(\nu)}/{\nu_j^{(k)}}+\mu_{j}^{(k)}}\tilde H_{i,j+1}^{(k)}(s)+\\
\nonumber
& \qquad \quad
\frac{\mu_{j}^{(k)}}{\lambda {f_{j}^{(k)}(\nu)}/{\nu_j^{(k)}}+\mu_{j}^{(k)}}\tilde H_{i-1,j-1}^{(k)}(s)\Bigg)\qquad  (2\leq i\leq j\leq B^{(k)}),\\
\label{eq:Hs}
\tilde H_{1,j}^{(k)}(s) &= \frac{\lambda f_{j}^{(k)}(\nu)/{\nu_j^{(k)}}+\mu_{j}^{(k)}}{s+\lambda f_{j}^{(k)}(\nu)/{\nu_j^{(k)}}+\mu_{j}^{(k)}}\Bigg(
\frac{\lambda {f_{j}^{(k)}(\nu)}/{\nu_j^{(k)}}}{\lambda {f_{j}^{(k)}(\nu)}/{\nu_j^{(k)}}+\mu_{j}^{(k)}}\tilde H_{1,j+1}^{(k)}(s)+\\
\nonumber
& \qquad \quad
\frac{\mu_{j}^{(k)}}{\lambda {f_{j}^{(k)}(\nu)}/{\nu_j^{(k)}}+\mu_{j}^{(k)}}
\Bigg) \qquad (1\leq j\leq B^{(k)}).
\end{align}
The corresponding version of \eqref{eq:Hweights} is 
\begin{align}
\label{eq:Hsweights}
\tilde H(s)= \sum_{k=1}^K\sum_{j=1}^{B^{(k)}} f^{(k)}_{j-1}(\nu)\tilde H_{j,j}^{(k)}(s).
\end{align}

Once again, \eqref{eq:Hs} and \eqref{eq:Hsweights} are valid when the functions $f_i^{(k)}$ are continuous at $\nu$. In other cases, we may need to tweak the formulas on a case-by-case basis.

The system time distribution can then be computed in the following manner:
\begin{enumerate}
\item We first compute the mean-field stationary distribution $\nu$. This can be done either by solving the balance equations \eqref{eq:dynamicbalance2}, or by numerically solving the mean-field transient equations \eqref{eq:mf_diff}, and setting a large enough $t$.
\item Once $\nu$ is available, \eqref{eq:Hs} is a system of linear equations for $\tilde H_{i,j}^{(k)}(s)$ that is straightforward to solve.
\item Then $\tilde H(s)$ is computed from \eqref{eq:Hsweights}.
\item Finally, $\tilde H(s)$ is transformed back to time domain.
\end{enumerate}
Due to \eqref{eq:Hs}, $\tilde H(s)$ is a rational function, whose inverse Laplace transform can be computed numerically. For numerical inverse Laplace transformation methods, we refer to \cite{ilt}.

We note that this approach to compute $\tilde H(s)$, while explicit, has its limitations, as the formula for $\tilde H(s)$ can get complicated for even moderately large $K$ and $B^{(k)}$ values. We address the feasibility further in Section \ref{ss:systemtimenum}.

Job losses occur only upon arrival, that is, all jobs that actually enter the system will be served, so $h_{i,j}^{(k)}(t)$ is a proper probability density function with
$$\int_0^\infty h_{i,j}^{(k)}(t) \mathrm{d}t=1.$$
However, if
$$\sum_{k=1}^K\sum_{j=1}^{B^{(k)}} f^{(k)}_{j-1}(\nu)<1,$$
then $\tilde H(s)$ is the Laplace-transform of a nonnegative function whose integral is equal to $1-\sum_{k=1}^K f^{(k)}_{B^{(k)}}(\nu)$ where
$$\sum_{k=1}^K f^{(k)}_{B^{(k)}}(\nu)$$
is the job loss probability, so in this sense, job losses are included in \eqref{eq:Hsweights}. The corresponding normalized version of \eqref{eq:Hsweights} is
\begin{align}
\label{eq:Hsnorm}
\frac{1}{1-\sum_{k=1}^K f^{(k)}_{B^{(k)}}(\nu)}\sum_{k=1}^K\sum_{j=1}^{B^{(k)}} f^{(k)}_{j-1}(\nu)\tilde H_{j,j}^{(k)}(s),
\end{align}
which is the Laplace-transform of a proper pdf whose integral is 1.

Depending on the load balancing principle, job losses may or may not be possible in the mean-field limit. This will be addressed specifically for each load balancing principle (For a finite system, job losses are always possible due to the finite buffers and fluctuations in either the job arrival or service speed.)

\section{Load balancing principles}
\label{s:loadbalance}
The load balancing principle describes the method the dispatcher uses to distribute the arriving jobs between the servers. It is quite important in large scale systems where the resources such as computing capacity are distributed between a large number of individual servers, and can make a big difference in the efficiency of the system.

The general goal of load balancing is to avoid long queues, directing incoming jobs to shorter queues instead.

There are several load balancing principles in use. Static policies do not consider the state of the system, only focusing on the incoming jobs. One example would be the round-robin load balancing policy, where incoming jobs are directed to the next server cyclically. Static load balancing principles are generally easy to operate, as they require minimal communication with the servers. Out of the principles observed in this paper, Random assignment falls into this category.

Dynamic principles, which take into account the current state of the system, can be more efficient. In real clusters, there is a trade-off: complicated policies require more communication and computation, generating a higher overhead communication cost, but provide better balancing. That said, in the mathematical framework we present, the cost of communication overhead is not modeled. Including the cost of overhead communication to provide an analytical framework for more realistic models is subject to further research.

In some systems it may be possible to reassign jobs that have been already assigned to new servers. It might also be possible that several servers ``team up'' to serve a single job. In our setting, we do not explore these options, and stick to a scenario where all jobs are assigned to a single server immediately upon arrival. On the other hand, in addition to the usual FIFO service principle, the framework does allow for limited processor sharing (LPS), where a single server can serve multiple jobs simultaneously.

In this paper we will examine 5 load balancing principles:
\begin{itemize}
\item Random assignment, where jobs are distributed randomly. With this principle, there is no actual load balancing. This principle will serve mostly as a baseline for comparison.
\item Join-Idle-Queue, where jobs are directed to idle queues if possible. A relatively recent idea \cite{jiq}, further explored in \cite{Mitzenmacher}.
\item Join-Shortest-Queue, where jobs are directed to the server with the fewest number of jobs waiting in queue. One of the earliest load balancing policies that has been widely used for decades \cite{jsq}. It provides very even balancing, but at the cost of high overhead communication, as the dispatcher needs to keep track of the queue length in every single server at all times.
\item Join-Shortest-Queue($d$), where jobs are directed to the server with the fewest number of jobs waiting in queue from among $d$ servers selected randomly. Also referred to as power-of-$d$, this is a version of JSQ that aims to reduce communication overhead at the cost of less strict balancing. It has been thoroughly explored, and has certain asymptotical optimality properties already for $d=2$ \cite{jsqd}.
\item Join-Below-Threshold, where jobs are directed to servers with a queue length below a prescribed threshold \cite{jbt}.
\end{itemize}

All of the above principles are based on natural intuitions that aim towards directing jobs to shorter queues, but they differ in the details and execution of doing so. In this section, we overview these load balancing principles from the literature. We present a high-level mathematical framework based on the Poisson representation of Section \ref{s:cluster} that is applicable to all of them, with the only difference being the $f_i^{(k)}(.)$ functions.

For each load balancing policy, we identify $f_i^{(k)}(.)$, then write the mean-field equations corresponding to \eqref{eq:mf}. We also identify the mean-field stationary distribution $\nu$ whenever available explicitly.

In case the $f_i^{(k)}(.)$ functions are discontinuous at $\nu$, we also rewrite the formulas \eqref{eq:Hij} and \eqref{eq:Hweights} so that they can be used to compute the mean system time, and rewrite the formulas \eqref{eq:Hs} and \eqref{eq:Hsweights} for system time distribution.

\subsection{Random assignment}

This is the most simple principle that we observe, and it does not lead to any balancing. With this setup the queues basically operate, and thus can be analyzed independently of each other. For random assignment,
$$f_i^{(k)}(x)=x^{(k)}_i,\quad k\in \{ 1,\dots,K\} ,$$
and accordingly, the mean-field equation is
\begin{align}
\label{eq:mf_random_ih}
\begin{split}
v_{i}^{(k)}(t)=&\int_0^t \lambda v^{(k)}_{i-1}(s)\d s -\int_0^t \lambda v^{(k)}_i(s)\d s\\
&\quad +\int_0^t \mu_{i+1} v^{(k)}_{i+1}(s)\d s -\int_0^t \mu_i v^{(k)}_{i}(s)\d s.
\end{split}
\end{align}
The mean-field balance equations, obtained from \eqref{eq:dynamicbalance2}, are
\begin{align}
\label{eq:bal_random2}
\mu_i^{(k)}\nu_i^{(k)}=\lambda \nu_{i-1}^{(k)}\qquad k\in \{ 1,\dots,K\} ,\qquad i\in \{ 1,\dots,B\} .
\end{align}

Solving \eqref{eq:bal_random2} gives the mean-field stationary distribution
$$\nu^{(k)}_i=c_k\prod_{j=1}^{i}\lambda/{\mu^{(k)}_j},\qquad i\in \{ 0,\dots,B^{(k)} \} ,$$
with the $c_k$'s coming from \eqref{eq:servertypes}. This is in accordance with the queues being independent.

Since the rates $f_i^{(k)}$ are continuous, \eqref{eq:Hij} and \eqref{eq:Hweights} can be used to compute the mean system time $H$, and \eqref{eq:Hs} and \eqref{eq:Hsweights} can be used to compute the Laplace-transform of the pdf of the system time distribution.

Job loss is possible for Random assignment, but is taken into account by the formulas \eqref{eq:Hs} and \eqref{eq:Hsweights}.

\subsection{Join-Idle-Queue}
\label{ss:jiq}

For Join-Idle-Queue (JIQ), incoming jobs are assigned to an idle server at random. If none of the servers are idle, a server is selected at random.

For JIQ, using the notation
$$y_0=\sum_{k=1}^K x_{0}^{(k)},$$
we have
\begin{align}
\label{eq:jiq_f}
f^{(k)}_i(x)=
\left\{
\begin{array}{ll}
\frac{x_{i}^{(k)}}{y_0} & \quad \textrm{ if }i=0,\,y_0 >0,\\
0 & \quad \textrm{ if }i>0,\,y_0 >0,\\
x_{i}^{(k)}&\quad \textrm{ if } y_0=0.
\end{array}
\right.
\end{align}
This system has been addressed in \cite{Mitzenmacher} for constant service rate curve and a homogeneous cluster.

The structure of the mean-field stationary distribution $\nu$ depends on the relation between $\lambda$ and $\sum_{k=1}^K \gamma_k \mu_1^{(k)}$. We address three cases separately.

\subsubsection*{JIQ, subcritical case}
When
$$\lambda<\sum_{k=1}^K \gamma_k \mu_1^{(k)},$$
there will always be idle queues in the mean-field stationary limit, so all jobs will be directed to idle queues. $\nu$ is concentrated on queues of length 0 and 1. From \eqref{eq:dynamicbalance2} we have
\begin{align}
\label{eq:bal_jiq}
\mu_1^{(k)} \nu_1^{(k)}=\lambda \frac{\nu_0^{(k)}}{\sum_{k=1}^K \nu_0^{(k)}}.
\end{align}

We do not have an explicit solution to \eqref{eq:bal_jiq}, but it can be solved numerically, and numerical experiments suggest a single fixed point $\nu$. In this region, the functions $f_i$ are continuous, so \eqref{eq:Hij} and \eqref{eq:Hweights} can be used to compute the mean system time $H$:
$$H=\sum_{k=1}^K \frac{\nu_0^{(k)}}{\sum_{k=1}^K \nu_0^{(k)}} H_{1,1}^{(k)},$$
and \eqref{eq:Hs} and \eqref{eq:Hsweights} can be used to compute the entire Laplace-transform of the system time distribution.

For subcritical JIQ, in the mean-field limit, there will be no job loss.

\subsubsection*{JIQ, critical case}
For
$$\lambda=\sum_{k=1}^K \gamma_k \mu_1^{(k)},$$ the mean-field stationary distribution is concentrated on queues of length 1, so we simply have
\begin{align}
\label{eq:jiq_crit}
\nu_1^{(k)}= \gamma_k, \quad k\in (1,\dots,K).
\end{align}

The functions $f_i^{(k)}$ are discontinuous at $\nu$, so \eqref{eq:Hij} and \eqref{eq:Hweights} does not apply. Instead, in the dynamic balance, whenever a queue of length 1 finishes service, a new job will enter immediately. With this, we can write the equivalent of \eqref{eq:Hij} for JIQ:
\begin{align}
\nonumber
H_{i,j}^{(k)} &= \frac{1}{\mu_{j}^{(k)}}+H_{i-1,j-1}^{(k)}
\qquad  (2\leq i\leq j\leq B^{(k)}),\\
\label{eq:Hij_jiq}
H_{1,j}^{(k)} &=\frac{1}{\mu_{j}^{(k)}} \qquad (1\leq j\leq B^{(k)}-1),
\end{align}

As we can see it is basically equivalent with \eqref{eq:Hij} in this case, because the discontinuity would only affect the arrival rate, and it is multiplied by 0 for every relevant term.
In the mean-field limit, all jobs go to queues of length 0 (which will then stay at length 1 for a positive amount of time), and there are no queues with 2 or more jobs. Accordingly, instead of \eqref{eq:Hweights}, we have
\begin{align}
\label{eq:H_jiq_crit}
H=\sum_{k=1}^K \frac{\mu_1^{(k)}\nu_1^{(k)}}{\lambda} H_{1,1}^{(k)}.
\end{align}

For the Laplace transforms, we have
\begin{align}
\nonumber
\tilde H_{i,j}^{(k)}(s) &= \frac{\mu_{j}^{(k)}}{s+\mu_{j}^{(k)}}
\tilde H_{i-1,j-1}^{(k)}(s),\qquad (2\leq i\leq j\leq B^{(k)}),\\
\label{eq:Hs_jiq_crit}
\tilde H_{1,j}^{(k)}(s) &=\frac{\mu_{j}^{(k)}}{s+\mu_{j}^{(k)}} \qquad (1\leq j\leq B^{(k)}-1),
\end{align}
and
\begin{align}
\label{eq:Hsweights_jiq_crit}
\tilde H(s)=\sum_{k=1}^K \frac{\mu_1^{(k)}\nu_1^{(k)}}{\lambda} \tilde H_{1,1}^{(k)}(s).
\end{align}

For critical JIQ, in the mean-field limit, there will be no job loss.

\subsubsection*{JIQ, supercritical case}
In case $\lambda>\sum_{k=1}^K \gamma_k \mu_1^{(k)}$, there will be no idle queues, so $\nu_0^{(k)}=0$ for $k\in (1,\dots,K)$. We note that $f_i^{(k)}$ are discontinuous at any point with $\sum_{k=1}^K\nu_0^{(k)}=0$ and $\sum_{k=1}^K\nu_1^{(k)}>0$; an intuitive explanation of this discontinuity is the following. Whenever a server with a single job finishes service, it will become idle. In the mean-field limit, a job will enter the idle queue instantly, so once again, we do not observe idle queues for any positive amount of time. However, similar to the $\lambda=\sum_{k=1}^K \gamma_k \mu_1^{(k)}$ case, a positive percentage of all incoming jobs will go to an idle queue. To compute this percentage, we once again observe that in the mean-field stationary distribution, service from queues of length 1 has to be balanced out completely by arrivals to idle queues.

The total service rate in queues of type $k$ of length 1 is $\mu_1^{(k)}\nu_1^{(k)}$, which is thus completely balanced out by an equal amount of arrivals The remaining arrival rate $(\lambda-\sum_{k=1}^K \mu_1^{(k)}\nu_1^{(k)})$ is distributed randomly. For longer queues, there are no discontinuities. Accordingly, the dynamic balance equations are
\begin{align}
\label{eq:jiq_mf_stat}
\left(\lambda-\sum_{k=1}^K \mu_1^{(k)}\nu_1^{(k)}\right)\nu_i^{(k)} = \mu_{i+1}^{(k)}\nu_{i+1}^{(k)},\quad i\in(1,\dots,B^{(k)}-1).
\end{align}

The system \eqref{eq:jiq_mf_stat} is nonlinear, but can be solved numerically. Then we can write a modified version of \eqref{eq:Hij} for the calculation of $H^{(k)}_{i,j}$. For this, we introduce
$$z_0=\sum_{k=1}^K \mu_1^{(k)}\nu_1^{(k)},$$
dubbed \emph{the upkeep}, which is the rate of service in servers with queue length 1, balanced out instantly by new arrivals. Essentially, the difference between \eqref{eq:jiq_mf_stat} and the original balance equations \eqref{eq:dynamicbalance2} is the presence of this upkeep term in the case when the dispatch functions are discontinuous at the mean-field stationary distribution $\nu$.

According to JIQ policy, the remaining arrival rate $\lambda-z_0$ is distributed randomly for the rest of the system. Accordingly, \eqref{eq:Hij} becomes
\begin{align}
\nonumber
H_{i,j}^{(k)} &= \frac{1}{(\lambda-z_0)+\mu_{j}^{(k)}}+
\frac{(\lambda-z_0)}{(\lambda-z_0)+\mu_{j}^{(k)}}H_{i,j+1}^{(k)}+\\
\nonumber
& \qquad \quad
\frac{\mu_{j}^{(k)}}{(\lambda-z_0)+\mu_{j}^{(k)}}H_{i-1,j-1}^{(k)}
\qquad  (2\leq i\leq j\leq B^{(k)}-1),\\
\label{eq:Hij_jiq_morethan1}
H_{i,B^{(k)}}^{(k)} &=\frac{1}{\mu_{B^{(k)}}^{(k)}}+H_{i-1,B^{(k)}-1}^{(k)} \qquad (2\leq i\leq B^{(k)}),\\
\nonumber
H_{1,j}^{(k)} &=\frac{1}{(\lambda-z_0)+\mu_{j}^{(k)}}+\frac{(\lambda-z_0)}{(\lambda-z_0)+\mu_{j}^{(k)}}H_{1,j+1}^{(k)} \quad (1\leq j\leq B^{(k)}-1),\\
\nonumber
H_{1,B^{(k)}}^{(k)} &=\frac{1}{\mu_{B^{(k)}}^{(k)}}.
\end{align}

To obtain the mean system time $H$, instead of \eqref{eq:Hweights}, we now have
\begin{align}
\label{eq:H_jiq}
H=\sum_{k=1}^K\frac{\mu_1^{(k)}\nu_1^{(k)}}{\lambda} H_{1,1}^{(k)}+ \left(1-\sum_{k=1}^K\frac{\mu_1^{(k)}\nu_1^{(k)}}{\lambda}\right)\sum_{k=1}^K\sum_{j=2}^{B^{(k)}}\nu_{j-1}^{(k)}H_{j,j}^{(k)}
\end{align}
since $\frac{\sum_{k=1}^K \mu_1^{(k)}\nu_1^{(k)}}{\lambda}$ is the portion of the arrival rate that is used to balance out the service in queues of length 1 and the remaining portion of the incoming rate is distributed randomly.

The corresponding equations for the Laplace transforms are
\begin{align}
\nonumber
\tilde H_{i,j}^{(k)}(s) &= \frac{(\lambda-z_0)+\mu_{j}^{(k)}}{s+(\lambda-z_0)+\mu_{j}^{(k)}}\Bigg(
\frac{(\lambda-z_0)}{(\lambda-z_0)+\mu_{j}^{(k)}}\tilde H_{i,j+1}^{(k)}(s)+\\
\nonumber
& \qquad \quad
\frac{\mu_{j}^{(k)}}{(\lambda-z_0)+\mu_{j}^{(k)}}\tilde H_{i-1,j-1}^{(k)}(s)\Bigg)
\qquad  (2\leq i\leq j\leq B^{(k)}-1),\\
\label{eq:Hsij_jiq_morethan1}
\tilde H_{i,B^{(k)}}^{(k)}(s) &=\frac{\mu_{B^{(k)}}^{(k)}}{s+\mu_{B^{(k)}}^{(k)}}\tilde H_{i-1,B^{(k)}-1}^{(k)}(s) \qquad (2\leq i\leq B^{(k)}),\\
\nonumber
\tilde H_{1,j}^{(k)}(s) &=\frac{(\lambda-z_0)+\mu_{j}^{(k)}}{s+(\lambda-z_0)+\mu_{j}^{(k)}}\Bigg(\frac{(\lambda-z_0)}{(\lambda-z_0)+\mu_{j}^{(k)}}\tilde H_{1,j+1}^{(k)}(s)+\\
\nonumber
& \qquad \quad \frac{\mu_{j}^{(k)}}{(\lambda-z_0)+\mu_{j}^{(k)}}\Bigg)\quad (1\leq j\leq B^{(k)}-1),\\
\nonumber
\tilde H_{1,B^{(k)}}^{(k)}(s) &=\frac{\mu_{B^{(k)}}^{(k)}}{s+\mu_{B^{(k)}}^{(k)}},
\end{align}
and
\begin{align}
\label{eq:Hsweights_jiq_supercrit}
\tilde H(s)=\sum_{k=1}^K\frac{\mu_1^{(k)}\nu_1^{(k)}}{\lambda} \tilde H_{1,1}^{(k)}(s)+ \left(1-\frac{z_0}{\lambda}\right)\sum_{k=1}^K\sum_{j=2}^{B^{(k)}}\nu_{j-1}^{(k)}\tilde H_{j,j}^{(k)}(s).
\end{align}

In general, for the supercritical JIQ case, job loss is possible, and is taken into account by the formula \eqref{eq:Hsweights_jiq_supercrit}.

\subsection{Join-Shortest-Queue}

For Join-Shortest-Queue (JSQ), incoming jobs are assigned to the shortest queue from among all queues; in case of multiple shortest queues of the same length, one is selected randomly.

For JSQ,
$$f^{(k)}_i(x)=
\left\{
\begin{array}{ll}
0 & \quad \textrm{ if }\exists\, i'<i\,\,\exists\, k':\, x_{i'}^{(k')}>0,\\
0 & \quad \textrm{ if }\sum_{k=1}^K x_{i}^{(k)}=0,\\
\frac{x_{i}^{(k)}}{\sum_{k=1}^K x_{i}^{(k)}}& \quad \textrm{ otherwise}.
\end{array}
\right.$$

For the stationary mean-field analysis, let $i_0$ denote the smallest $i$ for which
$$\sum_{k=1}^K\gamma_k\mu^{(k)}_i\geq \lambda.$$
Such an $i$ exists if the stability condition \eqref{eq:stable} holds. Then the mean-field stationary distribution $\nu$ will be concentrated on queues of length $i_0$ and $i_0-1$: starting from an arbitrary point, queues shorter than $i_0-1$ will receive the entire load of arrivals, which is larger than they can process, so these queues will ``fill up'' to level $i_0-1$, while queues longer than $i_0$ do not receive any load at all, so these queues will go down, until they reach level $i_0$.

The upkeep term is very similar to the JIQ case. The total service rate in queues of length $(i_0-1)$ is
$$z_0=\sum_{k=1}^{K}\mu^{(k)}_{i_0-1}\nu^{(k)}_{i_0-1},$$ which is completely balanced out by an equal amount of arrivals. In case $i_0=1$, $z_0=0$, so there is no upkeep, and all queues are of length 0 or 1; in this case, JSQ is equivalent to either subcritical or critical JIQ. When $i_0>1$, there is an actual upkeep. We assume $i_0>1$ for the rest of this section.

The remaining arrival rate $(\lambda-z_0)$ goes to queues of length $i_0-1$, with the queue type $k$ chosen at random with probabilities proportional to $\nu^{(k)}_{i_0-1}$. For each server type $k$, these arrivals are balanced out by the service in queues of type $k$ and length $i_0$, leading to the balance equations
\begin{align}
\label{eq:jsq_mf_stat2}
\mu^{(k)}_{i_0}\nu^{(k)}_{i_0} =\left(\lambda-z_0\right) \frac{\nu^{(k)}_{i_0-1}}{\sum_{k=1}^K \nu^{(k)}_{i_0-1}} \qquad k\in (1,\dots,K),
\end{align}
which, along with \eqref{eq:servertypes}, give a (nonlinear) system of equations for $\nu$, which can be solved numerically.

Whenever a server with queue length $i_0-1$ finishes service, it will become the single shortest queue and receives a new arrival instantly. Rate $(\lambda-z_0)$ remains for the rest of the system, which will be directed entirely to queues of length $i_0-1$. To ease notation, we also introduce
$$y_0=\sum_{k=1}^{K}\nu^{(k)}_{i_0-1}.$$
Then
\begin{align}
\nonumber
H_{i,j}^{(k)} &= H_{i,j+1}^{(k)}
\qquad  (1\leq i\leq j<i_0-1),\\
\nonumber
H_{1,i_0-1}^{(k)} &= \frac{1}{((\lambda-z_0)/y_0) + \mu_{i_0-1}^{(k)}} +\\
\nonumber
& \qquad\qquad
\frac{(\lambda-z_0)/y_0}{((\lambda-z_0)/y_0) + \mu_{i_0-1}^{(k)}} H_{1,i_0},\\
\nonumber
H_{i,i_0-1}^{(k)} &= \frac{1}{((\lambda-z_0)/y_0) + \mu_{i_0-1}^{(k)}} +\\
\label{eq:Hij_jsq}
 & \,\,\,\qquad
\frac{(\lambda-z_0)/y_0}{((\lambda-z_0)/y_0) + \mu_{i_0-1}^{(k)}} H_{i,i_0} +\\
\nonumber
& \qquad\qquad
\frac{\mu_{i_0-1}^{(k)}}{((\lambda-z_0)/y_0) + \mu_{i_0-1}^{(k)}} H_{i-1,i_0-2} \quad (2\leq i \leq i_0-1),\\
\nonumber
H_{1,j}^{(k)} &=\frac{1}{\mu_{j}^{(k)}} \qquad (i_0-1< j\leq B^{(k)}),\\
\nonumber
H_{i,j}^{(k)} &=\frac{1}{\mu_{j}^{(k)}}+H_{i-1,j-1} \quad (i_0-1<  j\leq B^{(k)}, \quad 1\leq i \leq j).
\end{align}
The first equation in (\ref{eq:Hij_jsq}) addresses the fact that if a server has fewer than $i_0-1$ jobs in it, it will immediately fill up to $i_0-1$ jobs. We also adjust the effective arrival rate to $\lambda -z_0 $, similarly to JIQ. If $i_0=1$, the $f_i^{(k)}$ are continuous at $\nu$, so we can use \eqref{eq:Hij} instead of \eqref{eq:Hij_jsq}. If $i_0=2$, there will of course not be any equation with the condition $(2\leq i \leq i_0-1)$.\\
If the functions $f_i^{(k)}$ are continuous at $\nu$, we can use \eqref{eq:Hweights} to calculate the mean system time. In case $i_0=1$, $\nu$ is in the inside of a continuous domain of the functions $f^{(k)}_i$, so this is the case, and \eqref{eq:Hweights} simplifies to
$$H=\sum_{k=1}^K \frac{\nu_{0}^{(k)}}{\sum_{k=1}^K \nu_{0}^{(k)}}  H^{(k)}_{1,1}.$$
On the other hand, if $i_0>1$, the functions $f_i$ are not continuous at $\nu$, and \eqref{eq:Hweights} is not applicable; instead, we have
$$H = \sum_{k=1}^K \frac{\mu^{(k)}_{i_0-1}\nu^{(k)}_{i_0-1}}{\lambda} H^{(k)}_{i_0-1,i_0-1}  +
\left(1-\frac{z_0}{\lambda}\right)
\sum_{k=1}^K  \frac{\nu^{(k)}_{i_0-1}}{\sum_{k=1}^K \nu^{(k)}_{i_0-1}}H^{(k)}_{i_0,i_0}  .$$

The corresponding equations for the Laplace transforms are
\begin{align}
\nonumber
\tilde H_{i,j}^{(k)}(s) &= \tilde H_{i,j+1}^{(k)}(s)
\qquad  (1\leq i\leq j<i_0-1),\\
\nonumber
\tilde H_{1,i_0-1}^{(k)}(s) &= \frac{(\lambda-z_0)/y_0 + \mu_{i_0-1}^{(k)}}{s+(\lambda-z_0)/y_0+ \mu_{i_0-1}^{(k)}} *\\
\nonumber
& \,\,\,\qquad
\Bigg(
\frac{\mu_{i_0-1}^{(k)}}{(\lambda-z_0)/y_0 + \mu_{i_0-1}^{(k)}}+
\frac{(\lambda-z_0)/y_0}{(\lambda-z_0)/y_0 + \mu_{i_0-1}^{(k)}} \tilde  H_{1,i_0}(s)\Bigg)\\
\label{eq:Hs_jsq}
\tilde  H_{i,i_0-1}^{(k)}(s) &= \frac{(\lambda-z_0)/y_0 + \mu_{i_0-1}^{(k)}}{s+(\lambda-z_0)/y_0 + \mu_{i_0-1}^{(k)}} *\\
\nonumber
& \,\,\,\qquad
\Bigg(\frac{(\lambda-z_0)/y_0}{(\lambda-z_0)/y_0 + \mu_{i_0-1}^{(k)}} \tilde  H_{i,i_0}(s) +\\
\nonumber
& \qquad\qquad
\frac{\mu_{i_0-1}^{(k)}}{(\lambda-z_0)/y_0 + \mu_{i_0-1}^{(k)}} \tilde H_{i-1,i_0-2}(s)\Bigg) \quad (2\leq i \leq i_0-1)\\
\nonumber
\tilde H_{1,j}^{(k)}(s) &=\frac{\mu_{j}^{(k)}}{s+\mu_{j}^{(k)}} \qquad (i_0-1< j\leq B^{(k)}),\\
\nonumber
\tilde H_{i,j}^{(k)}(s) &=\frac{\mu_{j}^{(k)}}{s+\mu_{j}^{(k)}}*\tilde H_{i-1,j-1}(s) \quad (i_0-1<  j\leq B^{(k)}, \quad 1\leq i \leq j),
\end{align}
and
$$\tilde H(s) = \sum_{k=1}^K \frac{\mu^{(k)}_{i_0-1}\nu^{(k)}_{i_0-1}}{\lambda} \tilde H^{(k)}_{i_0-1,i_0-1}(s)  +
\left(1-\frac{z_0}{\lambda}\right)
\sum_{k=1}^K  \frac{\nu^{(k)}_{i_0-1}}{\sum_{k=1}^K \nu^{(k)}_{i_0-1}}\tilde H^{(k)}_{i_0,i_0}(s).$$

Since $y_0$ and $z_0$ are straightforward to compute from $\nu$, \eqref{eq:Hs_jsq} is still a linear system of equations for $\tilde H_{i,j}^{(k)}(s)$, which is not any more difficult to solve than \eqref{eq:Hs}.

For JSQ, there is no job loss in the mean-field limit. (We emphasize that this is due to the stability condition \eqref{eq:stable}, which we assume in all cases.)

\subsection{Join-Shortest-Queue($d$)}

JSQ($d$) is a version of JSQ where the dispatcher first selects $d$ servers randomly, and dispatches the incoming job to the shortest from among the $d$ queues.\\
If we set $d=1$, we get Random assignment, and if we set $d=N$, we get JSQ. The $f_i^{(k)}$ functions are continuous for any finite $d$. Appendix \ref{app:jsqdconv} addresses the case $d\to\infty$.

For JSQ($d$), we introduce the auxiliary variables
$$y_i^{(k),N}=\sum_{j=i}^{B^{(k)}}x_{j}^{(k),N},\qquad z_i^N=\sum_{k=1}^K y_i^{(k),N},$$
and then inclusion-exclusion shows
\begin{align*}&f^{(k),N}_i(x^N)=
\frac{x_i^{(k),N}}{\sum_{k=1}^K x_i^{(k),N}}\times\\
&\quad\bigg[z_i^N\left(z_i^N-\frac1N\right)\dots\left(z_i^N-\frac{d-1}N\right)
-z_{i+1}^N\left(z_{i+1}^N-\frac1N\right)\dots\left(z_{i+1}^N-\frac{d-1}N\right)\bigg].
\end{align*}

The above version of $f^N_i(.)$ is $N$-dependent, but converges to
\begin{align*}f_i^{(k)}(x)=\frac{x_i^{(k)}}{\sum_{k=1}^K x_i^{(k)}}((z_i)^d-(z_{i+1})^d).
\end{align*}

Due to the dependency on $N$, we refer to \cite{Benait1}, where this type of dependence on $N$ is allowed. Also, both $f_i^{(k),N}$ and $f_i^{(k)}$ are continuous. Overall, the conclusions of Theorems \ref{t:kurtz-trans} and \ref{t:kurtz-stat} apply.

The mean-field balance equations are 
\begin{align}
\label{eq:bal_jsqd}
\frac{\lambda \nu_i^{(k)}}{\sum_{k=1}^K \nu_i^{(k)}}
\left(\left(\sum_{k=1}^K\sum_{j=i}^{B^{(k)}}\nu_{j}^{(k)}\right)^d-\left(\sum_{k=1}^K\sum_{j=i+1}^{B^{(k)}}\nu_{j}^{(k)}\right)^d\right)
=\mu_i^{(k)}\nu_i^{(k)}.
\end{align}

Since the rates $f_i^{(k)}$ are continuous, \eqref{eq:Hij} and \eqref{eq:Hweights} can be used to compute the mean system time $H$, and \eqref{eq:Hs} and \eqref{eq:Hsweights} can be used to compute the Laplace-transform of the pdf of the system time distribution.

Job loss is possible for JSQ($d$), but will be typically small enough to be negligible in practice.

\subsection{Join-Below-Threshold}

Join-Below-Threshold (JBT) sets a threshold $M_k$ which may depend on the server type $k$; servers of type $k$ with queue length $<M_k$ are considered available and servers of type $k$ with queue length $\geq M_k$ are full. Tasks will be dispatched to a random available servers. If there are no available servers, jobs will be dispatched at random among all servers.

JBT is commonly used in accordance with limited processor sharing (LPS) for servers which can serve multiple jobs simultaneously in an efficient manner. This is reflected in an increasing service rate curve $\mu_i^{(k)}$. If $\mu^{(k)}_i$ would start to decrease for large $i$, this is countered by setting the threshold $M_k$ at the maximum point. $M_k$ is referred to as the multi programming level (MPL), and is the number of jobs served simultaneously in a single server, while further jobs wait in queue. Overall, this setup ensures the service rate curve $\mu^{(k)}_i$ is increasing up to $M_k$ and constant for $M_k\leq i\leq B^{(k)}$.

If we set the threshold to 1, we get the JIQ principle, and if we set it to $B^{(k)}$, we get Random assignment.

We introduce the auxiliary variable
$$y= \sum_{k=1}^K\sum_{j=0}^{M_k-1}x^{(k)}_{j},$$
which is the ratio of available servers.
For JBT, 
$$f_i^{(k)}(x)=
\left\{
\begin{array}{ll}
0 & \quad \textrm{ if } y>0,\, i\geq M_k,\\
x^{(k)}_i/y& \quad \textrm{ if }y>0, \, i<M_k,\\
x^{(k)}_i & \quad \textrm{ if } y=0.
\end{array}
\right.$$

The mean-field balance equations are
\begin{align*}
\mu^{(k)}_i \nu^{(k)}_i =\frac{\lambda \nu_{i-1}^{(k)}}{y},\qquad i\in \{1,\dots,M_k-1\} ,\, \qquad k\in \{ 1,\dots,K \},
\end{align*}
with $\nu_i^{(k)}=0$ for $i>M_k$.

For a full, detailed mean-field analysis of JBT, we refer to \cite{jbt}. Apart from the stability condition \eqref{eq:stable} and monotonicity condition \eqref{eq:mon}, it is usually also assumed that
\begin{align}
\lambda<\sum_{k=1}^K \gamma_k\mu_{M_k},
\end{align}
which is a stability condition stronger than \eqref{eq:stable}, ensuring that the evolution of the transient mean-field limit eventually enters and then never leaves the region where no queues are longer than the threshold. On this domain, the functions $f_i^{(k)}$ are continuous, and the mean-field stationary solution $\nu$ is unique and also inside this domain. An efficient numerical method to compute $\nu$ is provided in \cite{jbt}.

As a side note, \cite{jbt} also shows examples where \eqref{eq:mon} does not hold, and there are multiple attractors in the mean-field system corresponding to quasi-stationary states of a system with a finite $N$, and mean-field convergence fails completely.

If \eqref{eq:stable} and \eqref{eq:mon} hold, \eqref{eq:Hij} and \eqref{eq:Hweights} can be used to compute the mean system time $H$, and \eqref{eq:Hs} and \eqref{eq:Hsweights} can be used to compute the Laplace-transform of the pdf of the system time distribution.

Job loss is not possible for JBT.

\section{Numerical experiments}
\label{s:numexp}

We conducted several numerical experiments. These are by no means exhaustive, but should nevertheless display some interesting properties and allow for some numerical comparison of the various load balancing methods.

For several parameter setups, we examined simulations for various choices of $N$, and also computed the mean-field limit ($N=\infty$). Simulations were done in Python and symbolic computations were done in Wolfram Mathematica. The codes for both are available at \cite{code}. For the symbolic calculations, numerical inverse Laplace transform was used, for which packages are available at \cite{iltorg}.

Section \ref{ss:trans} displays transient mean-field convergence as $N$ is increased. Also, as $t$ is increased, each system will converge to its stationary state.

Section \ref{ss:mean} compares the mean service times for both simulations and the mean-field settings.

Section \ref{ss:systemtimenum} addresses service time distributions.


\subsection{Homogeneous transient mean-field diagrams}
\label{ss:trans}

In this section, we plot the solutions of the mean-field equations as well as the corresponding $x_i^{(k),N}$ curves for systems with $N=1000$ and $N=10000$ servers, resulting from simulations. 

We will focus on homogeneous clusters with $K=1$ (also dropping $(k)$ from the notation). $B=B^{(k)}$, the maximal queue length will be set to 10. The rest of the parameter setup is shown in Table \ref{t:homogenparams}. The parameter setup adheres to the monotonicity assumption (\ref{eq:mon}) and also the stability condition (\ref{eq:stable}) (in fact, the system load can be computed as $\lambda/\mu_B$ in a homogeneous cluster).

Figures \ref{f:homrandom}--\ref{f:homjbt} display simulation results for the transient evolution of the homogeneous system using various load balancing policies. For each load balancing policy, two plots are included: the number of servers is $N=1000$ for the plot on the left and $N=10000$ for the plot on the right. Other system parameters are according to Table \ref{t:homogenparams}. All systems are initially empty. The x axis is time, and the jagged line graphs show the ratio of servers with queue length 0 to 10 respectively. These have some natural fluctuations. Also included are the transient mean-field limits, which are smooth curves.

\begin{table}
\begin{center}
\begin{tabular}{| c | c | c| c | c| c| c|} 
 \hline
 $\lambda$ & $\mu_1$ & $\mu_2$ & $\mu_3$ & $\mu_4$ & $\mu_5$ & $\mu_{6},\dots,\mu_{10}$\\
 \hline
 \hline
 $1.25$ & $1$ & $1.1$ & $1.2$ & $1.3$ & $1.4$ & $1.5$ \\
 \hline
\end{tabular}
\end{center}
\caption{Parameter setup for the homogeneous systems}
\label{t:homogenparams}
\end{table}

\subsubsection*{Random}

\begin{figure}
\centering
\begin{subfigure}{0.49\textwidth}
  \centering
  \includegraphics[width=1\linewidth]{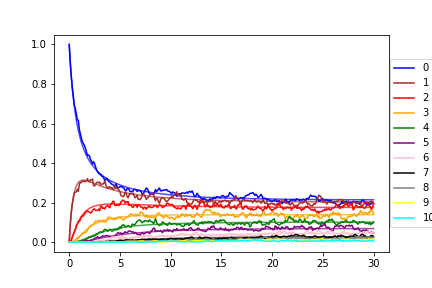}
  \caption{$N=1000$}
\end{subfigure}
\begin{subfigure}{0.49\textwidth}
  \centering
  \includegraphics[width=1\linewidth]{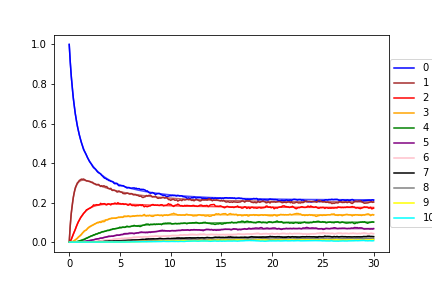}
  \caption{$N=10000$}
\end{subfigure}
\caption{Homogeneous transient evolution with Random load balancing}
\label{f:homrandom}
\end{figure}

Figure \ref{f:homrandom} displays the transient evolution with Random load balancing policy. A significant ratio of queues is longer throughout; overall, the Random load balancing principle is rather inefficient, and serves mostly as a baseline. Later we will see the effect of more efficient load balancing principles on the same systems.

The fluctuations of the simulations decrease as $N$ is increased. Actually, as mentioned after Theorem \ref{t:kurtz-trans}, the fluctuations are guaranteed to be of order $\frac{1}{\sqrt{N}}$ for $x^N$ (or, equivalently, order $\sqrt{N}$ for $X^N$). However, the constant factor can be different for the various load balancing principles. For Random assignment, the fluctuations are relatively mild.

Convergence to stationarity can also be observed: as time increases, the smooth graphs converge to the mean-field stationary distribution. That said, for any fixed finite $N$, the order of the fluctuations will not go to 0 as time is increased.

\subsubsection*{JIQ}

\begin{figure}[t]
\centering
\begin{subfigure}{0.49\textwidth}
  \centering
  \includegraphics[width=1\linewidth]{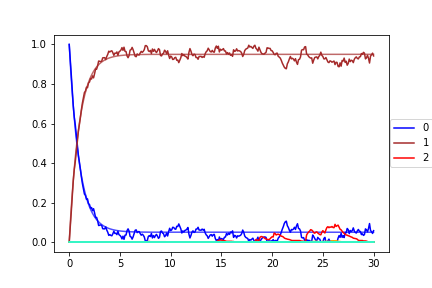}
  \caption{$N=1000,\lambda=0.95$}
  \label{f:homjiq1000_1}
\end{subfigure}
\begin{subfigure}{0.49\textwidth}
  \centering
  \includegraphics[width=1\linewidth]{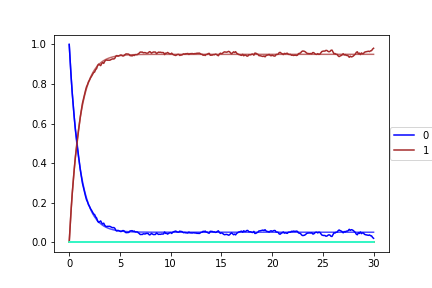}
  \caption{$N=10000,\lambda=0.95$}
  \label{f:homjiq10000_1}
\end{subfigure}
\\
\centering
\begin{subfigure}{0.49\textwidth}
  \centering
  \includegraphics[width=1\linewidth]{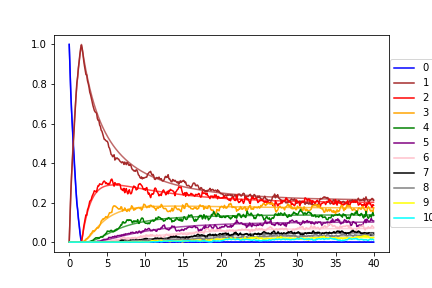}
  \caption{$N=1000,\lambda=1.25$}
  \label{f:homjiq1000_2}
\end{subfigure}
\begin{subfigure}{0.49\textwidth}
  \centering
  \includegraphics[width=1\linewidth]{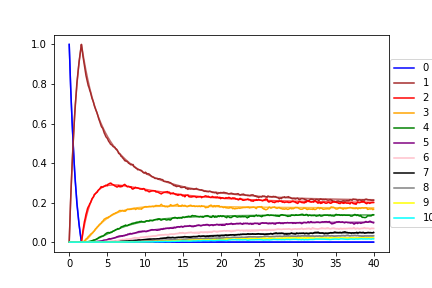}
  \caption{$N=10000,\lambda=1.25$}
  \label{f:homjiq10000_2}
\end{subfigure}
\caption{Homogeneous transient evolution with JIQ load balancing}
\label{f:homjiq}
\end{figure}

Figure \ref{f:homjiq} displays the transient evolution with JIQ load balancing policy for $\lambda=0.95$ and $\lambda=1.25$.

Figures \ref{f:homjiq1000_1} and \ref{f:homjiq10000_1} have $\lambda=0.95$ (with other parameters according to Table \ref{t:homogenparams}), which is subcritical due to $\lambda=0.95<\mu_1=1$ (see Section \ref{ss:jiq}), so the system stabilizes on queues of length 0 and 1.

Figures \ref{f:homjiq1000_2} and \ref{f:homjiq10000_2} have $\lambda=1.25>\mu_1=1$, which is supercritical, so the system starts out by filling up all empty queues in a sharp manner. After this initial period, no empty queues are present anymore, and the dynamic dispatch is distributed among queues of length 1 through 10 randomly. Similar to Random policy, once again longer queues are present in the system. 

\subsubsection*{JSQ(2) and JSQ(5)}

\begin{figure}
\centering
\begin{subfigure}{0.49\textwidth}
  \centering
  \includegraphics[width=1\linewidth]{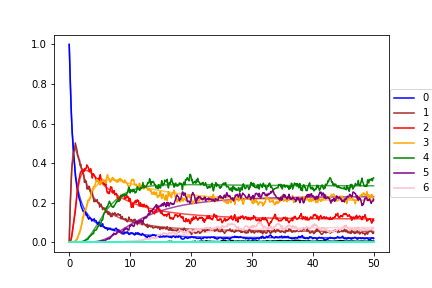}
  \caption{$N=1000$}
\end{subfigure}
\begin{subfigure}{0.49\textwidth}
  \centering
  \includegraphics[width=1\linewidth]{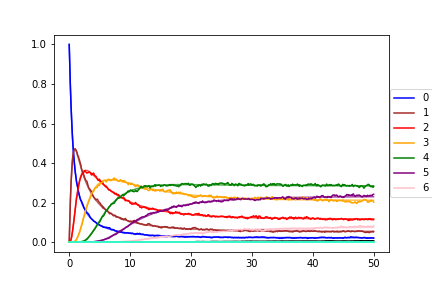}
  \caption{$N=10000$}
\end{subfigure}
\caption{Homogeneous transient evolution with JSQ(2) load balancing}
\label{f:homjsq2}
\end{figure}

Figure \ref{f:homjsq2} displays the transient evolution with JSQ(2) load balancing policy. Already for $d=2$, the result is markedly different from Random assignment. This is a known phenomenon, referred to as power-of-2 \cite{jsqd}. The ratio of longer queues diminishes more rapidly with the queue length than for either Random or JIQ policy.

\begin{figure}
\centering
\begin{subfigure}{0.49\textwidth}
  \centering
  \includegraphics[width=1\linewidth]{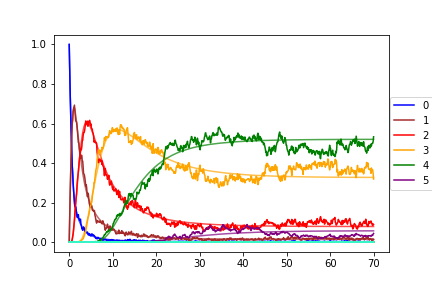}
  \caption{$N=1000$}
\end{subfigure}
\begin{subfigure}{0.49\textwidth}
  \centering
  \includegraphics[width=1\linewidth]{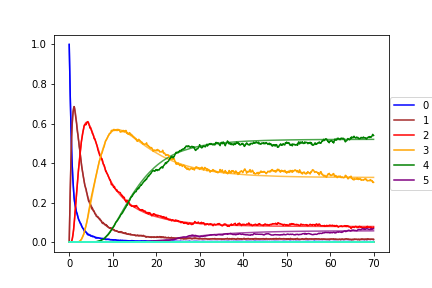}
  \caption{$N=10000$}
\end{subfigure}
\caption{Homogeneous transient evolution with JSQ(5) load balancing}
\label{f:homjsq5}
\end{figure}

Figure \ref{f:homjsq5} displays the transient evolution with JSQ(5) load balancing policy. Here, most of the queues will be of length 3 and 4, with the ratio of either shorter or longer queues much smaller. We also note that the dispatch function is continuous, so the transient mean-field limit functions are smooth, although they change rather sharply.

\subsubsection*{JSQ}

\begin{figure}
\centering
\begin{subfigure}{0.49\textwidth}
  \centering
  \includegraphics[width=1\linewidth]{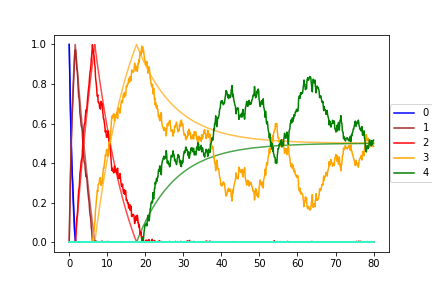}
  \caption{$N=1000$}
\end{subfigure}
\begin{subfigure}{0.49\textwidth}
  \centering
  \includegraphics[width=1\linewidth]{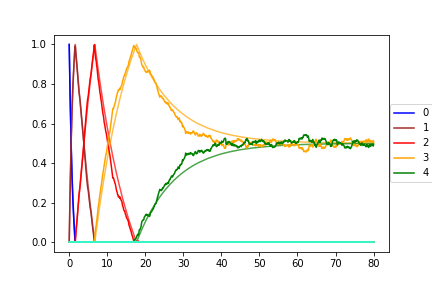}
  \caption{$N=10000$}
\end{subfigure}
\caption{Homogeneous transient evolution with JSQ load balancing}
\label{f:homjsq}
\end{figure}

Figure \ref{f:homjsq} displays the transient evolution with JSQ load balancing policy. Here, all of the queues will be of length 3 and 4 after the system fills up. At any point in time, there are only 2 different queue lengths present, starting from lengths 0 and 1, switching to 1 and 2, then 2 and 3, then 3 and 4 as the system fills up. We also note that the dispatch function is discontinuous, so the transient mean-field limit functions has breaking points at switches to new queue length pairs.

The stationary mean-field limit is $\nu_3=\nu_4=0.5$ due to $$\lambda=1.25=\frac{\mu_3+\mu_4}{2}=\frac{1.2+1.3}{2}.$$

For any finite $N$, when a job in a queue of minimal length finishes service, a shorter queue will appear for a brief but positive time. In the mean-field limit, such queues are filled back instantly.

We also note that the fluctuations are considerably larger than for either Random or JIQ. An intuitive explanation is that the higher level of control provided by JSQ will generally focus any fluctuations in either the arrival or service on a single queue length: if the arrivals outweigh the service for a short period of time, the surplus arrivals will all go to servers of minimal queue length. Overall, the strict control introduces a positive correlation between the length of the queues, resulting in larger fluctuations (which are, once again, of order $1/\sqrt{N}$, but with a higher constant factor). Principles with less strict control generally distribute this fluctuation among several different queue lengths, resulting in smaller fluctuations.

\subsubsection*{JBT}

\begin{figure}
\centering
\begin{subfigure}{0.49\textwidth}
  \centering
  \includegraphics[width=1\linewidth]{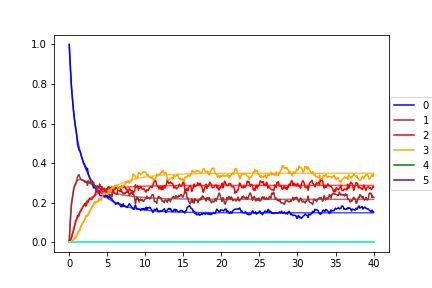}
 \caption{$N=1000$}
\end{subfigure}
\begin{subfigure}{0.49\textwidth}
  \centering
  \includegraphics[width=1\linewidth]{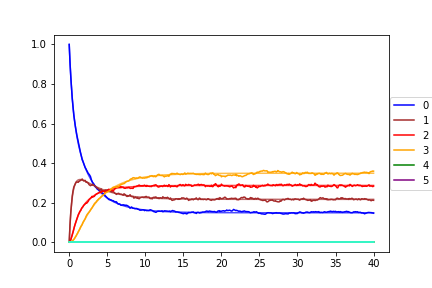}
 \caption{$N=10000$}
\end{subfigure}
\caption{Homogeneous transient evolution with JBT load balancing}
\label{f:homjbt}
\end{figure}

Figure \ref{f:homjbt} displays the transient evolution with JBT load balancing policy. The MPL parameter is set to 5. In this setup, the system reaches stability before hitting the MPL threshold (and accordingly, the mean-field system reaches its attractor before the discontinuity point, so the functions remain continuous). This is the intended usage of JBT.

\subsection{Heterogeneous transient mean-field diagrams}
\label{ss:heterogentrans}

In this section, we plot the solutions of the mean-field equations as well as the corresponding $x_i^{(k),N}$ curves for systems with $N=10000$ servers, resulting from simulations. 

We will focus on heterogeneous clusters with $K=2$. $B=B^{(k)}$, the maximal queue length will be set to 10. The rest of the parameter setup is shown in Table \ref{t:heterogenparams}. The parameter setup adheres to the monotonicity assumption (\ref{eq:mon}) and also the stability condition (\ref{eq:stable}).

The parameter choices in Table \ref{t:heterogenparams} are motivated by an actual real-life scenario: in many shopping centers, there are two types of checkouts: checkouts served by an employee (service rate $1$ in Table \ref{t:heterogenparams}), with a separate queue for each such checkout, and self-service checkouts. A single self-service checkout is typically slightly slower (service rate $0.8$ in Table \ref{t:heterogenparams}) than a checkout served by an employee, but this is countered by the fact that there is a batch of self-service checkouts for each queue (the batch size is 5 for Table \ref{t:heterogenparams}).

Of course, in actual shopping centers, the number of queues may or may not be high enough to warrant a mean-field approach; that said, as we will see later, some derived performance measures are well-approximated by the mean-field limit already for smaller system sizes.

Figures \ref{f:hetrandom}--\ref{f:hetjbt} display simulation results for the transient evolution of the heterogeneous system using various load balancing policies. For each load balancing policy, two plots are included: the ratio of type $1$ servers with various queue lengths for the plot on the left and the ratio of type $2$ servers with various queue lengths for the plot on the right. Other system parameters are according to Table \ref{t:heterogenparams}. All systems are initially empty. The x axis is time, and the jagged line graphs show the ratio of servers with queue length 0 to 10 respectively. These have some natural fluctuations. Also included are the transient mean-field limits, which are smooth curves.

\begin{table}[t]
\begin{center}
\begin{tabular}{|c || c|c | c| c | c| c| c|} 
 \hline
 $k$ & $\lambda$ & $\mu_1^{(k)}$ & $\mu_2^{(k)}$ & $\mu_3^{(k)}$ & $\mu_4^{(k)}$ & $\mu_5^{(k)}$ & $\mu_{6}^{(k)},\dots,\mu_{10}^{(k)}$\\
 \hline
 \hline
 $1$ & \multirow{2}{*}{$1.6$} & $1.0$ & $1.0$ & $1.0$ & $1.0$ & $1.0$ & $1.0$ \\ 
 \cline{1-1} \cline{3-8}
 $2$ & & $0.8$ & $1.6$ & $2.4$ & $3.2$ & $4.0$ & $4.0$ \\ 
 \hline
\end{tabular}
\end{center}
\caption{Parameter setup for the heterogeneous systems}
\label{t:heterogenparams}
\end{table}

\subsubsection*{Random}

\begin{figure}[ht]
\centering
\begin{subfigure}{0.49\textwidth}
  \centering
  \includegraphics[width=1\linewidth]{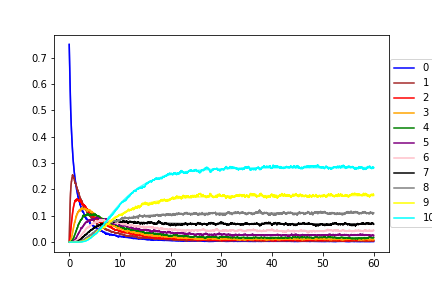}
  \caption{type 1}
\end{subfigure}
\begin{subfigure}{0.49\textwidth}
  \centering
  \includegraphics[width=1\linewidth]{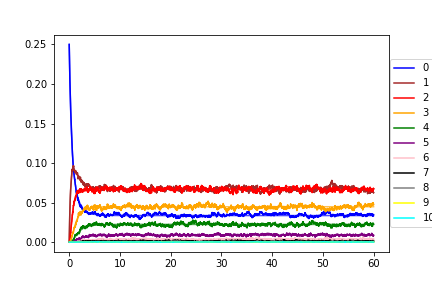}
  \caption{type 2}
\end{subfigure}
\caption{Heterogeneous transient evolution with Random load balancing}
\label{f:hetrandom}
\end{figure}

Figure \ref{f:hetrandom} displays the transient evolution with Random load balancing policy. A significant ratio of queues is longer throughout; in fact, servers of type 1 are overloaded, as can be seen from the fact that the majority of queues of type 1 has length 10 (equal to the buffer size) or close. In a heterogeneous system, with poor load balancing, it is possible that some server types are overloaded even though the system as a whole is subcritical.

\subsubsection*{JIQ}

\begin{figure}[ht]
\centering
\begin{subfigure}{0.49\textwidth}
  \centering
  \includegraphics[width=1\linewidth]{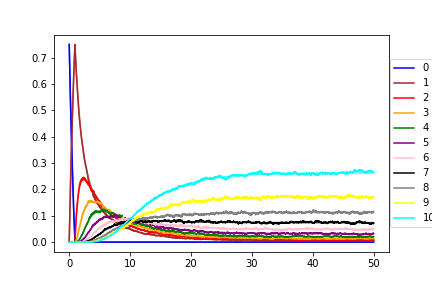}
  \caption{server type 1}
\end{subfigure}
\begin{subfigure}{0.49\textwidth}
  \centering
  \includegraphics[width=1\linewidth]{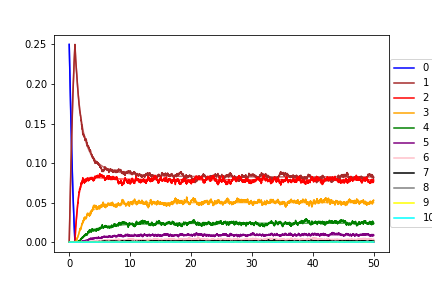}
  \caption{server type 2}
\end{subfigure}
\caption{Heterogeneous transient evolution with JIQ load balancing}
\label{f:hetjiq}
\end{figure}

Figure \ref{f:hetjiq} displays the transient evolution with JIQ load balancing policy.

JIQ does not offer a considerable improvement over Random, as once again longer queues are present in the system. This also means that servers of type 1 are overloaded, which also results in significant data loss. On the other hand, servers of type 2 are subcritical.

\subsubsection*{JSQ($2$) and JSQ($5$)}

\begin{figure}[ht]
\centering
\begin{subfigure}{0.49\textwidth}
  \centering
  \includegraphics[width=1\linewidth]{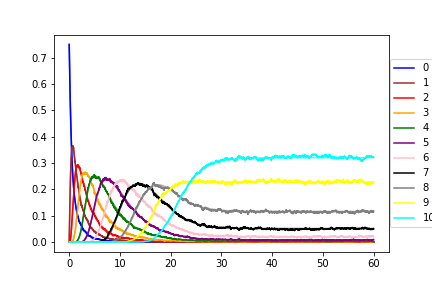}
  \caption{server type 1}
\end{subfigure}
\begin{subfigure}{0.49\textwidth}
  \centering
  \includegraphics[width=1\linewidth]{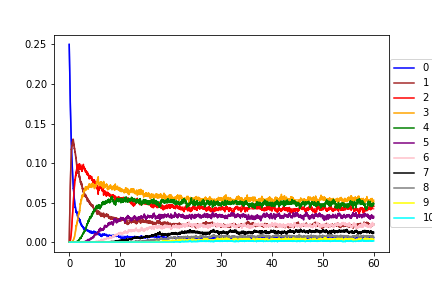}
  \caption{server type 2}
\end{subfigure}
\caption{Heterogeneous transient evolution with JSQ($2$) load balancing}
\label{f:hetjsq2}
\end{figure}

Figure \ref{f:hetjsq2} displays the transient evolution with JSQ($2$) load balancing policy. Servers of type 1 are still overloaded, in which case JSQ(2) does not offer a considerable improvement over either Random or JIQ. The system (particularly servers of type 1) goes through an initial build-up period, starting from empty and converging to stationarity with the majority of queues full (length equal to buffer size 10) or close.

\begin{figure}[ht]
\centering
\begin{subfigure}{0.49\textwidth}
  \centering
  \includegraphics[width=1\linewidth]{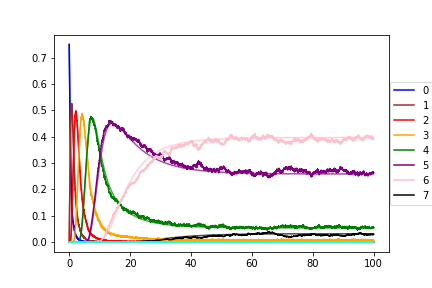}
  \caption{server type 1}
\end{subfigure}
\begin{subfigure}{0.49\textwidth}
  \centering
  \includegraphics[width=1\linewidth]{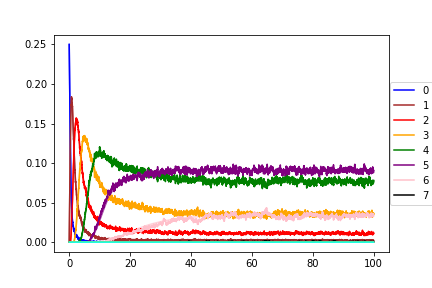}
  \caption{server type 2}
\end{subfigure}
\caption{Heterogeneous transient evolution with JSQ($5$) load balancing}
\label{f:hetjsq5}
\end{figure}

Figure \ref{f:hetjsq5} displays the transient evolution with JSQ($5$) load balancing policy. In this case, the better load balancing results in both server types being subcritical; for server type 1, the typical queue lengths are 5 and 6, while for server type 2, the typical queue lengths are 4 and 5. Data loss is practically negligible in this case.

\subsubsection*{JSQ}

\begin{figure}[ht]
\centering
\begin{subfigure}{0.49\textwidth}
  \centering
  \includegraphics[width=1\linewidth]{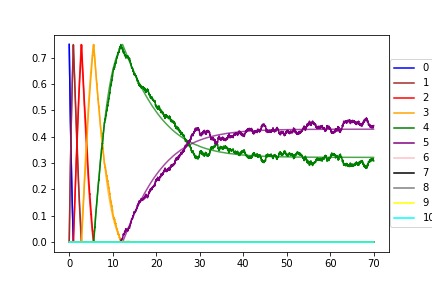}
  \caption{server type 1}
\end{subfigure}
\begin{subfigure}{0.49\textwidth}
  \centering
  \includegraphics[width=1\linewidth]{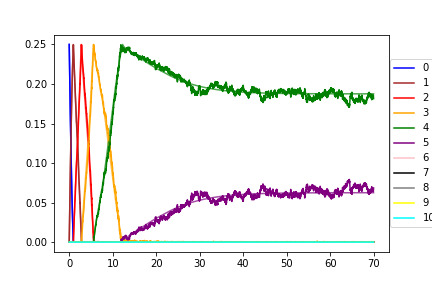}
  \caption{server type 2}
\end{subfigure}
\caption{Heterogeneous transient evolution with JSQ load balancing}
\label{f:hetjsq}
\end{figure}

Figure \ref{f:hetjsq} displays the transient evolution with JSQ load balancing policy. The build-up period is much sharper (in fact, the mean-field limit curves are nondifferentiable at the changes in minimal queue length), with both server types eventually reaching a state where all queue lengths are either 4 or 5. Fluctuations around the mean-field limit are relatively mild for $N=10000$ servers.

\subsubsection*{JBT}

\begin{figure}[ht]
\centering
\begin{subfigure}{0.49\textwidth}
  \centering
  \includegraphics[width=1\linewidth]{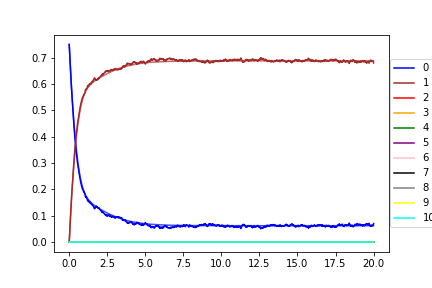}
  \caption{server type 1}
\end{subfigure}
\begin{subfigure}{0.49\textwidth}
  \centering
  \includegraphics[width=1\linewidth]{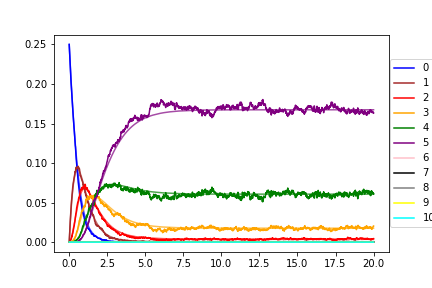}
  \caption{server type 2}
\end{subfigure}
\caption{Heterogeneous transient evolution with JBT load balancing}
\label{f:hetjbt}
\end{figure}

Figure \ref{f:hetjbt} displays the transient evolution with JBT load balancing policy. MPL parameters are 1 for server type 1 and 5 for server type 2. JBT load balancing policy suits the type of heterogeneous system described by Table \ref{t:heterogenparams} particularly well: the MPL settings allow to fully utilize the service capacity of each server type without allowing queues longer than necessary. In fact, JBT can outperform JSQ for heterogeneous systems, as we will see in the next section.

\subsection{Mean system times}
\label{ss:mean}

The main performance measure we are going to examine is the mean system time, that is, the average time a job spends between arrival and finishing service.

First we examine the homogeneous system described by the parameter settings in Table \ref{t:homogenparams} for simulations for various system sizes ranging from $N=10$ to $N=10000$ and also the mean-field limit, with the various load balancing principles from Section \ref{s:loadbalance}. Table \ref{t:hommean} lists the mean system times from both simulations, and calculated from the mean-field limit using equations (\ref{eq:Hij}) and (\ref{eq:Hweights}) (or in the discontinuous cases, their corresponding versions listed in Section \ref{s:loadbalance}). We note that despite long running times, the simulation results still may have an inherent small random variation.

\begin{table}[ht]
	\begin{center}
\begin{tabular}{|c||c|c|c|c|c|}
\hline
{Load balancing} & $N=10$ & $N=100$ & $N=1000$ & $N=10000$ & $N=\infty$ \\ \hline
Random  & 3.560 & 3.563 & 3.571 & 3.560 & 3.565 \\ \hline
JIQ     & 3.034 & 2.914 & 2.907 & 2.903 & 2.886 \\ \hline
JSQ(2)  & 3.031 & 2.963 & 2.961 & 2.960 & 2.958 \\ \hline
JSQ(5)  & 2.857 & 2.813 & 2.819 & 2.817 & 2.817 \\ \hline
JSQ     & 2.834 & 2.792 & 2.802 & 2.800 & 2.800 \\ \hline 
JBT     & 3.087 & 3.010 & 2.996 & 2.995 & 2.993 \\ \hline
\end{tabular}
\end{center}
\caption{Mean system time in the stationary mean-field limit (homogeneous cluster)}
\label{t:hommean}
\end{table}

JSQ is the most effective principle, which is unsurprising (although we do emphasize that in practice, JSQ comes with a heavy overhead communication burden which was not modelled here). 

JSQ($d$) is more effective with a higher $d$, but already for $d=2$, it is significantly better than Random, which is once again known as the power-of-$2$ (or power-of-$d$) \cite{mitzenmacher2}.

We note that jobs lost are not included in the averages in Table \ref{t:hommean}; in order to give a more complete picture, we mention that the theoretical job loss probability for Random policy (with the same parameters as per Table \ref{t:homogenparams}) is $0.0438$, and for JIQ it is $0.0136$ (for JSQ(2), JSQ(5), JSQ and JBT, job loss is negligible). Job loss probabilities for the simulations are not included in the paper, we just mention that they closely match the theoretical values.

Overall, based on Table \ref{t:hommean}, the mean-field approximation for the mean system times is exceedingly accurate already for small values of $N$.

Next we address the heterogeneous system described by the parameter settings in Table \ref{t:heterogenparams}.

As long as $N$ is finite, there are fluctuations which do not vanish even as time increases and the systems converge to their stationary limit. As expected, fluctuations are bigger for smaller values of $N$. For smaller values of $N$, the mean system time is generally above the mean-field mean system time; an intuitive explanation for this is that the limited number of servers offers less `room' to balance out short periods of overflow (coming from the natural fluctuations of arrivals and service), causing the system to operate with longer queues for said short periods.

\begin{table}[ht]
	\begin{center}
\begin{tabular}{|c||c|c|c|c|c|c|}
\hline
\multirow{2}{*}{Load bal.}  & \multirow{2}{*}{Server type} & \multicolumn{5}{c|}{$N$} \\
\cline{3-7} 
& & $12$ & $100$ \! &\! $1000$ & $10000$ & $\infty$ \\ \hline
\hline
\multirow{3}{*}{Random}
& Entire system  & 5.935 & 5.934      & 5.936       & 5.932         & 5.933 \\ \cline{2-7} 
& Server type 0  & 8.434 & 8.427      & 8.432       & 8.423         & 8.425 \\ \cline{2-7} 
& Server type 1  & 1.274 & 1.274      & 1.273       & 1.274        & 1.274 \\  \hline \hline
\multirow{3}{*}{JIQ}
& Entire system  & 5.651 & 5.631      & 5.639       & 5.640         & 5.638 \\ \cline{2-7} 
& Server type 0  & 8.265 & 8.233      & 8.244       & 8.251         & 8.246 \\ \cline{2-7} 
& Server type 1  & 1.272 & 1.268      & 1.269       & 1.271         & 1.270 \\ \hline \hline
\multirow{3}{*}{JSQ(2)} 
& Entire system  & 5.217 & 5.347      & 5.353       & 5.348         & 5.352 \\ \cline{2-7} 
& Server type 0  & 8.794 & 8.979      & 8.982       & 8.972         & 8.976 \\ \cline{2-7} 
& Server type 1  & 1.408 & 1.389      & 1.382       & 1.383         & 1.381 \\  \hline \hline
\multirow{3}{*}{JSQ(5)}
& Entire system  & 3.390 & 3.287      & 3.273       & 3.272         & 3.273 \\ \cline{2-7} 
& Server type 0  & 5.685 & 5.536      & 5.514       & 5.519         & 5.517 \\ \cline{2-7} 
& Server type 1  & 1.378 & 1.302      & 1.294       & 1.293         & 1.293 \\  \hline \hline
\multirow{3}{*}{JSQ}
& Entire system  & 3.086 & 2.797      & 2.806       & 2.807         & 2.807 \\ \cline{2-7} 
& Server type 0  & 5.084 & 4.556      & 4.569       & 4.573         & 4.571 \\ \cline{2-7} 
& Server type 1  & 1.336 & 1.249      & 1.250       & 1.250        & 1.250 \\  \hline \hline
\multirow{3}{*}{JBT}
& Entire system  & 2.568 & 1.175      & 1.142       & 1.143        & 1.143 \\ \cline{2-7} 
& Server type 0  & 1.298 & 1.071      & 0.999       & 1.001        & 1.000 \\ \cline{2-7} 
& Server type 1  & 1.856 & 1.253      & 1.250       & 1.250        & 1.250 \\  \hline                                
\end{tabular}
\end{center}
\caption{Mean system time in the stationary mean-field limit (heterogeneous cluster)}
\label{t:mean}
\end{table}

Once again, in order to compare the mean system time for the various load balancing principles, it is important to take into account that some of these principles operate with significant data loss: for random, the theoretical job loss probability is $0.285$, for JIQ, it is $0.251$, and for JSQ(2), it is $0.104$.

Table \ref{t:mean} shows that, similar to the homogeneous case (Table \ref{t:hommean}), the mean-field approximation for the mean system times is very accurate for both smaller and larger choices of $N$ (and for JSQ(5), JSQ and JBT, job loss is negligible). The only exception is JBT for $N=12$; for very small system sizes and system load close to critical ($1.6/1.75$ according to the parameters in Table \ref{t:heterogenparams}), even a small burst in the arrivals can push the entire system over the threshold, at which point it switches to Random, and stays there for significant periods of time.

\subsection{System time distributions}
\label{ss:systemtimenum}

In this section we examine the theoretical probability density function of the system time in the mean-field limit for some setups and compare it with empirical distributions (histograms) from simulations for finite $N$.

The theoretical distributions are calculated using equations \eqref{eq:Hs} and \eqref{eq:Hsweights} (or in discontinuous cases their counterparts described in Section \ref{s:loadbalance}), and inverse Laplace transformation (ILT). The system \eqref{eq:Hs} can be solved explicitly, and the solution is a rational function (in the Laplace transform domain).

However, depending on the value of $K$ and $B^{(1)},\dots,B^{(K)}$, the solution for $\tilde H(s)$ from \eqref{eq:Hsweights} can be infeasible already for moderately large values of $K$ and $B$. In general, the formula for $\tilde H(s)$ is relatively simple if only few of the $\tilde H_{i,j}^{(k)}$'s are nonzero, which is typically the case for JSQ. For other load balancing principles, where all $\tilde H_{i,j}^{(k)}$'s are nonzero, the explicit formula for $\tilde H(s)$ from \eqref{eq:Hsweights} is infeasible already for $K=2$ and $B^{(1)}=B^{(2)}=10$.

Due to this, the parameters for this setup were the homogeneous system from Table \ref{t:homogenparams} with $\lambda=1.25$. We also set $B=5$, to make the ILT less complicated. Just as an example, for JSQ, with the above parameters, we have
$$\tilde H(s)=\frac{(24 s+65)^4}{5 (2 s+5)^3 (10 s+13)^4}.$$
$\tilde H(s)$ can be computed for the other load balancing principles as well, but the explicit formulas are far more complicated, and are omitted from the paper.
\begin{figure}[h!t]
\centering
\begin{subfigure}{0.49\textwidth}
  \centering
  \includegraphics[width=1\linewidth]{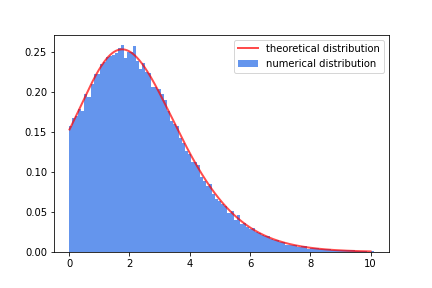}
  \caption{Random}
\end{subfigure}
\begin{subfigure}{0.49\textwidth}
  \centering
  \includegraphics[width=1\linewidth]{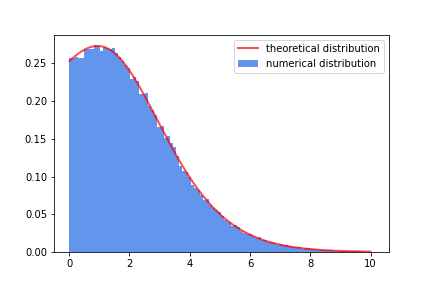}
  \caption{JIQ}
\end{subfigure}

\begin{subfigure}{0.49\textwidth}
  \centering
  \includegraphics[width=1\linewidth]{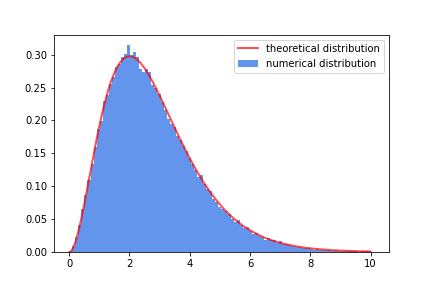}
  \caption{JSQ}
\end{subfigure}
\begin{subfigure}{0.49\textwidth}
  \centering
  \includegraphics[width=1\linewidth]{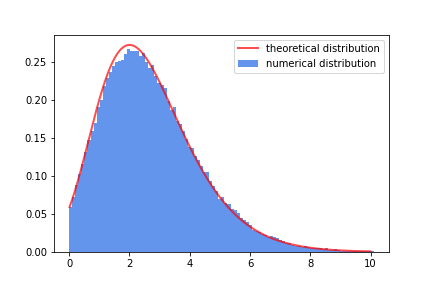}
  \caption{JSQ(2)}
\end{subfigure}

\begin{subfigure}{0.49\textwidth}
  \centering
  \includegraphics[width=1\linewidth]{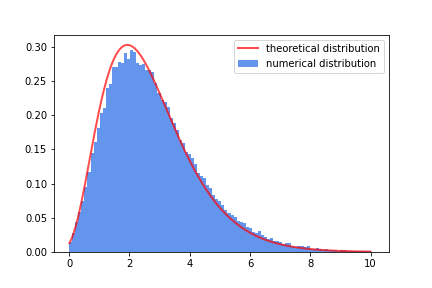}
  \caption{JSQ(5)}
\end{subfigure}
\begin{subfigure}{0.49\textwidth}
  \centering
  \includegraphics[width=1\linewidth]{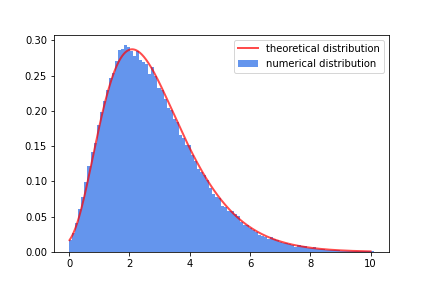}
  \caption{JBT}
\end{subfigure}
\caption{System time distributions}
\label{f:systemtime}
\end{figure}

Figure \ref{f:systemtime} displays the theoretical pdf of the system time in the mean-field limit with a red curve, while the blue histograms are from simulations with $N=1000$ servers. Each system was run long enough to reach the stationary regime, and only jobs arriving during this period were considered. The theoretical pdf's are normalized as per \eqref{eq:Hsnorm}.

In general, all histograms match the theoretical pdf's well. For random assignment and JIQ (which is supercritical with the given parameters), the system time is less concentrated (e.g. it has a higher variance). JSQ is the only one where the system time density is 0 at time $t=0$; for all other load balancing principles, it is possible that a job starts service immediately, which corresponds to a positive density at $t=0$. For JSQ($2$) and JSQ($5$), the match between the theoretical and numerical distributions is slightly less perfect than for others (although still very good); the exact reason for this is subject to further research.

\section{Conclusion and outlook}
\label{s:concl}

In this paper we examined the mean-field transient and stationary convergence of systems with several different load-balancing principles based on queue length. 

While no rigorous proof was presented, the simulations suggest that mean-field convergence holds even for discontinuous $f_i^{(k)}$ dispatch functions. We have provided formulas to compute the stationary mean-field limit, and also the mean system time in the mean-field stationary regime. In addition to that, the entire service time distribution could also be calculated with the help of the Laplace transform, adapting \eqref{eq:Hij} and \eqref{eq:Hweights} for the Laplace transforms of the system times. We have also examined the mean system time numerically for several parameter setups.

There is a lot of possibility for further work in this topic. One direction would be to provide mathematically rigorous proofs for versions of Theorems \ref{t:kurtz-trans} and \ref{t:kurtz-stat} for some of the discussed systems with discontinuous dispatch functions.

Another direction is scenarios where further information is available (e.g. job size); in such cases, that information can be used to estimate the load of each queue more precisely and design other load balancing principles. 

Yet another direction is to add a geometrical dimension to the server cluster, with the load balancing principle taking into account the distance of the arriving job to the queues (e.g. as in a shopping center, where customers are more likely to choose a queue physically closer to their arrival point).

We could also make the model more realistic, even if more complicated, by considering the dispatcher's communication overhead cost. However, we expect the communication overhead cost to be highly dependent on actual system settings, and as such, it seems difficult to incorporate it in a high level model in a general manner.

Another direction is to allow different job types, with certain job types can be served more efficiently by certain server types.

All in all, this is a vast topic that has a lot of potential for further development.

\bibliographystyle{abbrv}
\bibliography{lbmf}

\appendix

\section{Little's law}
\label{app:little}

In a heterogeneous system, Little's law applies to the entire system in the mean-field stationary regime, and also applies to each server type separately. It is valid regardless if the dispatch functions are continuous or not, but requires some consideration for discontinuous dispatch functions. In this section, we provide the proper formulas for each load balancing principle.

Let $\lambda^{(k)}$ denote the effective arrival rate to servers of type $k$, and $L^{(k)}$ denote the average queue length in servers of type $k$ ($k=1,\dots, K$). Using these, we can compute the mean system time for a job in a server of type $k$ via Little's law as $$H^{(k)}=L^{(k)}/\lambda^{(k)}.$$

For any load balancing principle, 
$$L^{(k)}=\frac{\sum_{i=0}^{B^{(k)}} i\nu_{i}^{(k)}}{\sum_{i=0}^{B^{(k)}} \nu_{i}^{(k)}}.$$

The formula for $\lambda^{(k)}$ is different for continuous and discontinuous dispatch functions. For dispatch functions continuous at $\nu$ (this case includes Random, JSQ($d$), JBT and also subcritical JIQ and JSQ with $i_0=1$), the formula for $\lambda^{(k)}$ is
$$\lambda^{(k)}=\lambda\frac{
\sum_{i=0}^{B^{(k)}-1}f_i^{(k)}(\nu)}{\sum_{i=0}^{B^{(k)}}\nu_i^{(k)}}.
$$
For supercritical JIQ, we have
$$\lambda^{(k)}=\frac{\mu_1^{(k)}\nu_1^{(k)}+(\lambda-z_0)
\sum_{i=1}^{B^{(k)}-1}\nu_i^{(k)}}{\sum_{i=0}^{B^{(k)}}\nu_i^{(k)}},
$$
and for JSQ with $i_0>1$, we have
$$\lambda^{(k)}=\frac{\mu_{i_0-1}^{(k)}\nu_{i_0-1}^{(k)}+(\lambda-z_0)
\frac{\nu_{i_0-1}^{(k)}}{\sum_{k=1}^K\nu_{i_0-1}^{(k)}}}{\nu_{i_0-1}^{(k)}+\nu_{i_0}^{(k)}}.
$$

\section{System time distribution for LPS service principle}
\label{app:lps}

This section is a counterpart of Section \ref{ss:systemtime}; we provide formulas to compute the system time distribution for limited processor sharing (LPS) service principle.

For LPS, each server type has a parameter called the \emph{multi-programming level} (MPL); the server can serve a number of jobs up to the MPL simultaneously, dividing its service capacity evenly, while further jobs wait in a FIFO queue.

Once again, let $h_{i,j}^{(k)}(t)$ denote the probability density function of the remaining system time of a job at position $i$ in a queue of length $j$ and queue type $k$. $M^{(k)}$ denotes the multi-programming level of queues of type $k$. The order of jobs is irrelevant among jobs already in service; that is, for fixed $k$ and $j$, $h_{i,j}^{(k)}(t)$ is constant for $i\leq \min(j,M^{(k)}).$ Accordingly, in the formulas we will write $h_{1,j}^{(k)}(t)$ instead of $h_{i,j}^{(k)}(t)$ for $i\leq \min(j,M^{(k)}).$ For jobs that are not yet in service ($i> M^{(k)}$), their position within the queue is still relevant.

For LPS, when the tagged job is in service, three type of changes can occur to its queue: arrival, or the tagged job finishes service, or another job finishes service. In the last case, it does not matter whether the finished job is ahead or behind the tagged job. When the tagged job is not yet in service, only two type of changes can occur: arrival, or another job finishes service. We also use once again that arrival is not possible when the queue is full ($j=B^{(k)}$), that is, $f_{B^{(k)}}^{(k)}(.)\equiv 0$ for $k=1,\dots,K$.

The corresponding version of \eqref{eq:Hs} is as follows:
\begin{align}
\nonumber
\tilde H_{1,j}^{(k)}(s) &= \frac{\lambda f_{j}^{(k)}(\nu)/{\nu_j^{(k)}}+\mu_{j}^{(k)}}{s+\lambda f_{j}^{(k)}(\nu)/{\nu_j^{(k)}}+\mu_{j}^{(k)}}\Bigg(
\frac{\lambda {f_{j}^{(k)}(\nu)}/{\nu_j^{(k)}}}{\lambda {f_{j}^{(k)}(\nu)}/{\nu_j^{(k)}}+\mu_{j}^{(k)}}\tilde H_{1,j+1}^{(k)}(s)+\\
\nonumber
& \qquad \quad
\frac{\mu_{j}^{(k)}(M^{(k)}-1)/M^{(k)}}{\lambda {f_{j}^{(k)}(\nu)}/{\nu_j^{(k)}}+\mu_{j}^{(k)}}\tilde H_{1,j-1}^{(k)}(s)+
\frac{\mu_{j}^{(k)}/M^{(k)}}{\lambda {f_{j}^{(k)}(\nu)}/{\nu_j^{(k)}}+\mu_{j}^{(k)}}\Bigg)\\
\label{eq:Hslps}
& \quad \qquad (1\leq i\leq M^{(k)}\leq j\leq B^{(k)}),\\
\nonumber
\tilde H_{1,j}^{(k)}(s) &= \frac{\lambda f_{j}^{(k)}(\nu)/{\nu_j^{(k)}}+\mu_{j}^{(k)}}{s+\lambda f_{j}^{(k)}(\nu)/{\nu_j^{(k)}}+\mu_{j}^{(k)}}\Bigg(
\frac{\lambda {f_{j}^{(k)}(\nu)}/{\nu_j^{(k)}}}{\lambda {f_{j}^{(k)}(\nu)}/{\nu_j^{(k)}}+\mu_{j}^{(k)}}\tilde H_{1,j+1}^{(k)}(s)+\\
\nonumber
& \qquad \quad
\frac{\mu_{j}^{(k)}(j-1)/j}{\lambda {f_{j}^{(k)}(\nu)}/{\nu_j^{(k)}}+\mu_{j}^{(k)}}\tilde H_{1,j-1}^{(k)}(s)+
\frac{\mu_{j}^{(k)}/j}{\lambda {f_{j}^{(k)}(\nu)}/{\nu_j^{(k)}}+\mu_{j}^{(k)}}\Bigg)\\
\nonumber
& \quad \qquad (1\leq i \leq j< M^{(k)}),
\end{align}
\begin{align}
\nonumber
\tilde H_{M^{(k)}+1,j}^{(k)}(s) &= \frac{\lambda f_{j}^{(k)}(\nu)/{\nu_j^{(k)}}+\mu_{j}^{(k)}}{s+\lambda f_{j}^{(k)}(\nu)/{\nu_j^{(k)}}+\mu_{j}^{(k)}}\Bigg(
\frac{\lambda {f_{j}^{(k)}(\nu)}/{\nu_j^{(k)}}}{\lambda {f_{j}^{(k)}(\nu)}/{\nu_j^{(k)}}+\mu_{j}^{(k)}}\tilde H_{M^{(k)}+1,j+1}^{(k)}(s)+\\
\nonumber
& \qquad \quad
\frac{\mu_{j}^{(k)}}{\lambda {f_{j}^{(k)}(\nu)}/{\nu_j^{(k)}}+\mu_{j}^{(k)}}\tilde H_{1,j-1}^{(k)}(s)\Bigg)
\quad ( j\leq B^{(k)}),\\
\label{eq:Hslps2}
\tilde H_{i,j}^{(k)}(s) &= \frac{\lambda f_{j}^{(k)}(\nu)/{\nu_j^{(k)}}+\mu_{j}^{(k)}}{s+\lambda f_{j}^{(k)}(\nu)/{\nu_j^{(k)}}+\mu_{j}^{(k)}}\Bigg(
\frac{\lambda {f_{j}^{(k)}(\nu)}/{\nu_j^{(k)}}}{\lambda {f_{j}^{(k)}(\nu)}/{\nu_j^{(k)}}+\mu_{j}^{(k)}}\tilde H_{i,j+1}^{(k)}(s)+\\
\nonumber
& \qquad \quad
\frac{\mu_{j}^{(k)}}{\lambda {f_{j}^{(k)}(\nu)}/{\nu_j^{(k)}}+\mu_{j}^{(k)}}\tilde H_{i-1,j-1}^{(k)}(s)\Bigg)\\
\nonumber
& \quad \qquad (M^{(k)}+1<i\leq j\leq B^{(k)}).
\end{align}

Once again, \eqref{eq:Hslps} and \eqref{eq:Hsweights} are applicable to compute $\tilde H(s)$ when the dispatch functions $f_i^{(k)}$ are continuous at $\nu$. In other cases, the formulas may need to be modified. 

\section{Partial control}
\label{app:partial}

We highlight a situation dubbed \emph{partial control}. In such a system, some of the jobs are not subject to the load balancing policy, and will simply be dispatched randomly. A real life example for partial control would be directing traffic via cooperating navigation apps in cars: each car with a cooperating navigation app is subject to load balancing, but drivers without the app select routes not subject to the same load balancing.

Assume we have a system with a load balancing policy corresponding to some dispatch functions $f_i^{(k)}(x)$. Load balancing only has partial control: for each job, with some fixed probability $0<p\leq 1$, the job will be dispatched according to the load balancing policy, but with probability $(1-p)$, it will be dispatched randomly. In this case, the corresponding dispatch functions are simply
$$\hat f_i^{(k)}(x) = p f_i^{(k)}(x) + (1-p)x_i^{(k)}.$$

\begin{figure}
\centering
\begin{subfigure}{0.49\textwidth}
  \centering
  \includegraphics[width=1\linewidth]{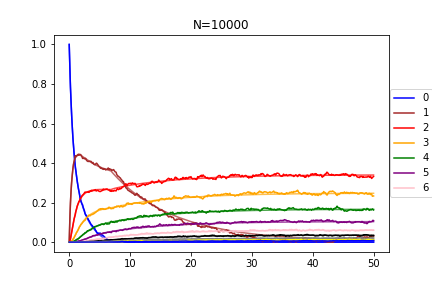}
  \caption{$p=0.3$}
  \label{f:jsq03}
\end{subfigure}
\begin{subfigure}{0.49\textwidth}
  \centering
  \includegraphics[width=1\linewidth]{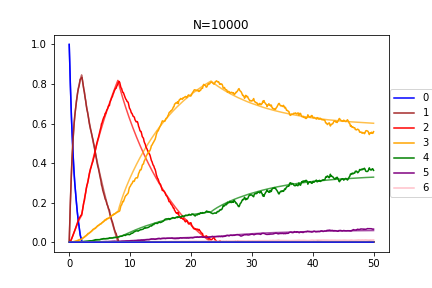}
  \caption{$p=0.8$}
  \label{f:jsq08}
\end{subfigure}
\caption{Partially controlled JSQ}
\label{fig:jsqpartial}
\end{figure}

Figure \ref{fig:jsqpartial} shows transient plots with JSQ load balancing principle with low ($p=0.3$) and high ($p=0.8$) levels of control. System parameters are according to Table \ref{t:homogenparams} with $\lambda = 1.25$ and $N=10000$. With a low level of control, the transient behaviour is closer to the case of random assignment, with longer queues also present. For low control, the minimal stationary queue length is 2, lower than the minimal stationary queue length 3 in case of full control JSQ, as the system needs to balance fewer controlled jobs (e.g.\ the upkeep is lower). For high control ($p=0.8$), the minimal stationary queue length remains 3, but once again, longer queues are also present.

\section{Convergence of JSQ($d$) to JSQ as $d\to\infty$}
\label{app:jsqdconv}

\begin{figure}
\centering
\begin{subfigure}{0.49\textwidth}
  \centering
  \includegraphics[width=1\linewidth]{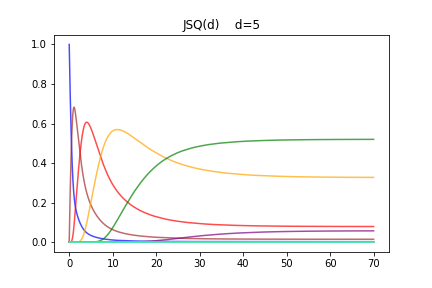}
\end{subfigure}
\begin{subfigure}{0.49\textwidth}
  \centering
  \includegraphics[width=1\linewidth]{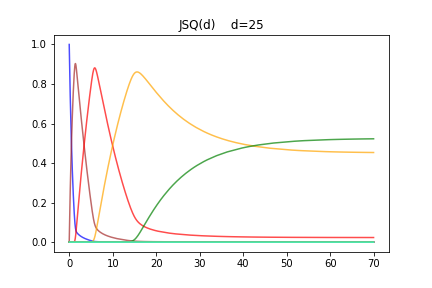}
\end{subfigure}
\begin{subfigure}{0.49\textwidth}
  \centering
  \includegraphics[width=1\linewidth]{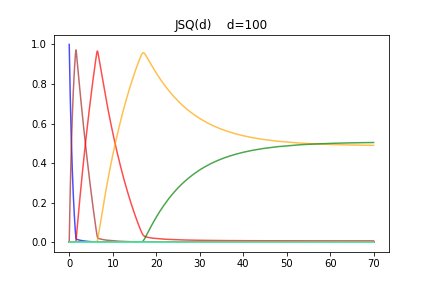}
\end{subfigure}
\begin{subfigure}{0.49\textwidth}
  \centering
  \includegraphics[width=1\linewidth]{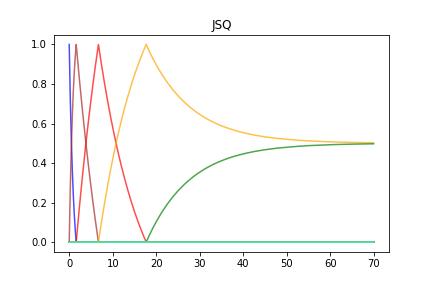}
\end{subfigure}
\caption{JSQ($d$)'s convergence to JSQ}
\label{fig:jsqdconv}
\end{figure}

This section shows an interesting visualisation of JSQ($d$)'s ``convergence'' to JSQ as $d\rightarrow \infty$. Figure \ref{fig:jsqdconv} displays the solutions of the transient mean-field equations for various choices of $d$. In practice, JSQ($d$) is quite close to JSQ already for moderately large values of $d$.

We note that the mean-field transient solutions are smooth for JSQ($d$) for any choice of $d$, but not for JSQ.

\end{document}